\documentclass[11 pt,Rectoverso]{article}  
\def\proof{\noindent{\bf Proof. }}

\def \1{\textrm{\dsrom{1}}}

\usepackage[latin1]{inputenc}
\usepackage{graphicx,psfrag,amsfonts,bbm}
\usepackage{amsmath,indentfirst,qsymbols}

\def \G{{\bf G}}
\def \oG{\overline{{\bf G}}}
\def \A{{\bf A}}

\def \be{\begin{eqnarray*}}
\def \ee{\end{eqnarray*}}
\def \ben{\begin{eqnarray}}
\def \een{\end{eqnarray}}
\def \bq{\begin{equation}}
\def \eq{\end{equation}}

\def \build#1#2#3{\mathrel{\mathop{\kern 0pt#1}\limits_{#2}^{#3}}}

\def \eref#1{(\ref{#1})}
\def \sous#1#2{\mathrel{\mathop{\kern 0pt#1}\limits_{#2}}}
\def \sur#1#2{\mathrel{\mathop{\kern 0pt#1}\limits^{#2}}}

\def \bG{{\bf G}}

\def \Proof {{\bf Proof }}
\def \captionn#1{\begin{center}\begin{minipage}{15cm}\sf\caption{\small #1}\end{minipage}\end{center}}

\def \dis{\displaystyle}
\def \tend{\longrightarrow}
\def \Sq{{\sf Sq}}
\def \Cy{{\sf Cy}}

\def \b{\big}
\def \bar{\overline}

\def \l{\left}
\def \r{\right}

\newcommand{\bt}{{\bf t}}

\font\dsrom=dsrom10 scaled 1400
\newqsymbol{`I}{\textrm{\dsrom{1}}}
\newqsymbol{`P}{\mathbb{P}}
\newqsymbol{`Q}{\mathbb{Q}}
\newqsymbol{`R}{\mathbb{R}}
\newqsymbol{`Z}{\mathbb{Z}}
\newqsymbol{`E}{\mathbb{E}}
\newqsymbol{`e}{\varepsilon}
\newqsymbol{`a}{\alpha}
\newqsymbol{`t}{\tau}
\newqsymbol{`o}{\omega}
\newqsymbol{`O}{\Omega}
\newqsymbol{`N}{\mathbb{N}}

\graphicspath{{../DESSINS/}}

\begin{document}
\newtheorem{lem}{Lemma}[section]
\newtheorem{defi}[lem]{Definition}
\newtheorem{pro}[lem]{Proposition}
\newtheorem{theo}[lem]{Theorem}
\newtheorem{cor}[lem]{Corollary}
\newtheorem{remi}[lem]{Remark\rm}{\rm}
\newtheorem{com}[lem]{Comments\rm}{\rm}
\newtheorem{exe}[lem]{Examples \rm}{\rm}

\begin{center}\LARGE{\bf Directed animals in the gas}\medskip\normalsize
\[\begin{array}{ll}
\textrm{\Large Yvan Le Borgne}& \textrm{\Large ~~~~~~~~~~Jean-Fran\c{c}ois Marckert}\end{array}\]
\textrm{CNRS, LaBRI}\\
\textrm{Universit\'e Bordeaux 1}\\
   \textrm{351 cours de la Libération}\\
 \textrm{33405 Talence cedex, France}
 \end{center}
\begin{abstract}
In this paper, we revisit the enumeration of directed animals using gas models. We show that there exists a natural construction of random directed animals on any directed graph together with a particle system that explains at the level of objects the formal link known between the density of the gas model and the generating function of directed animals counted according to the area. This provides some new methods to compute the generating function of directed animals counted according to area, and leads in the particular case of the square lattice to new combinatorial results and questions. A model of gas related to directed animals counted according to area and perimeter on any directed graph is also exhibited.
\end{abstract}

\section{Introduction}

\subsection{Directed animals on a directed graph}\label{gc}
Let $G=(V,E)$ be a connected graph with set of vertices $V$ and set of edges $E$. An animal $A$ in $G$ is a subset of $V$ such that between any two vertices $u$ and $v$ in $A$, there is a path in $G$ having all its vertices in $A$. The vertices of $A$ are called cells and the number of cells is denoted $|A|$ and called the area of $A$. A neighbor $u$ of $A$ is a vertex of $G$ which is not in $A$ and such that there exists $e\in E$ between $u$ and a vertex $v(u)$ of $A$. The perimeter of $A$, denoted by ${\cal P}(A)$, is the set of neighbors of $A$ and its cardinality is denoted by $|{\cal P}(A)|$.

Directed animals (DA) are animals built on a directed graph  $G$~: 
\begin{defi}\label{diran}
Let $A$ and $S$ be two subsets of $V$, with {\bf finite or infinite} cardinalities. We say that $A$ is a DA  with source $S$, if $S$ is a subset of $A$ such that any vertex of $A$ can be reached from an element of $S$ through a directed path having all its vertices in $A$.
\end{defi}
In the setting of DA, the definition of cells and area are the same as in the case of animals, but the notion of neighbor is changed since the edge $(v(u),u)$ is required to be a directed edge of $G$. If there is a directed edge from $v$ to $u$ in $G$, then $u$ is said to be a child of $v$, and $v$ to be the father of $u$; this induces a notion of descendant and ancestor. Each node of ${\cal P}(A)$ has at least one father in $A$. 
In this paper, we deal only with DA built on graphs having some suitable properties~:
\begin{defi} A directed graph $G=(V,E)$ is said to be agreeable if\\
(A) $G$ does not contain multiple edges,\\
(B) the graph $G$ has no directed cycles,\\
(C) the number of children of each node is finite.
\end{defi}
Condition (A) is needed to avoid some useless discussions in the Theorems. Condition (C) is needed to have a finite number of DA with a given area, for all sources. 
Notice that an agreeable graph is not necessarily connected, and is not necessarily locally finite (some nodes may have an infinite incoming degree). Even if never recalled in the statements of the results, $V$ is supposed to be finite or countable. \par
As examples, (finite or infinite) trees and forests are agreeable graphs, the square lattice  $\Sq=\mathbb{Z}^2$ directed in such a way that the vertex $(x,y)$ has as children $(x,y+1)$ and $(x+1,y+1)$ is agreeable (as well as all usual directed lattices).\medskip

 \noindent A subset $S$ of $V$ is said to be free if for any $x,y\in S$, $x\neq y$, $x$ is not an ancestor of $y$.
For any DA $A$, the set 
\[{\cal S}(A):=\{x, x\in A, x\textrm{ has no father in  }A\}\] 
is a free subset of $V$ and is the unique minimal source of $A$ according to the inclusion partial order (it is also the intersection of all possibles sources of $A$). \par

\begin{figure}[htbp]
\psfrag{A}{$a$}\psfrag{B}{$b$}\psfrag{C}{$c$}
\psfrag{G_0}{$G_0$}
\psfrag{G_1}{$G_1$}
\psfrag{G_2}{$G_2$}
\psfrag{G_3}{$G_3$}
\centerline{\includegraphics[height=6cm]{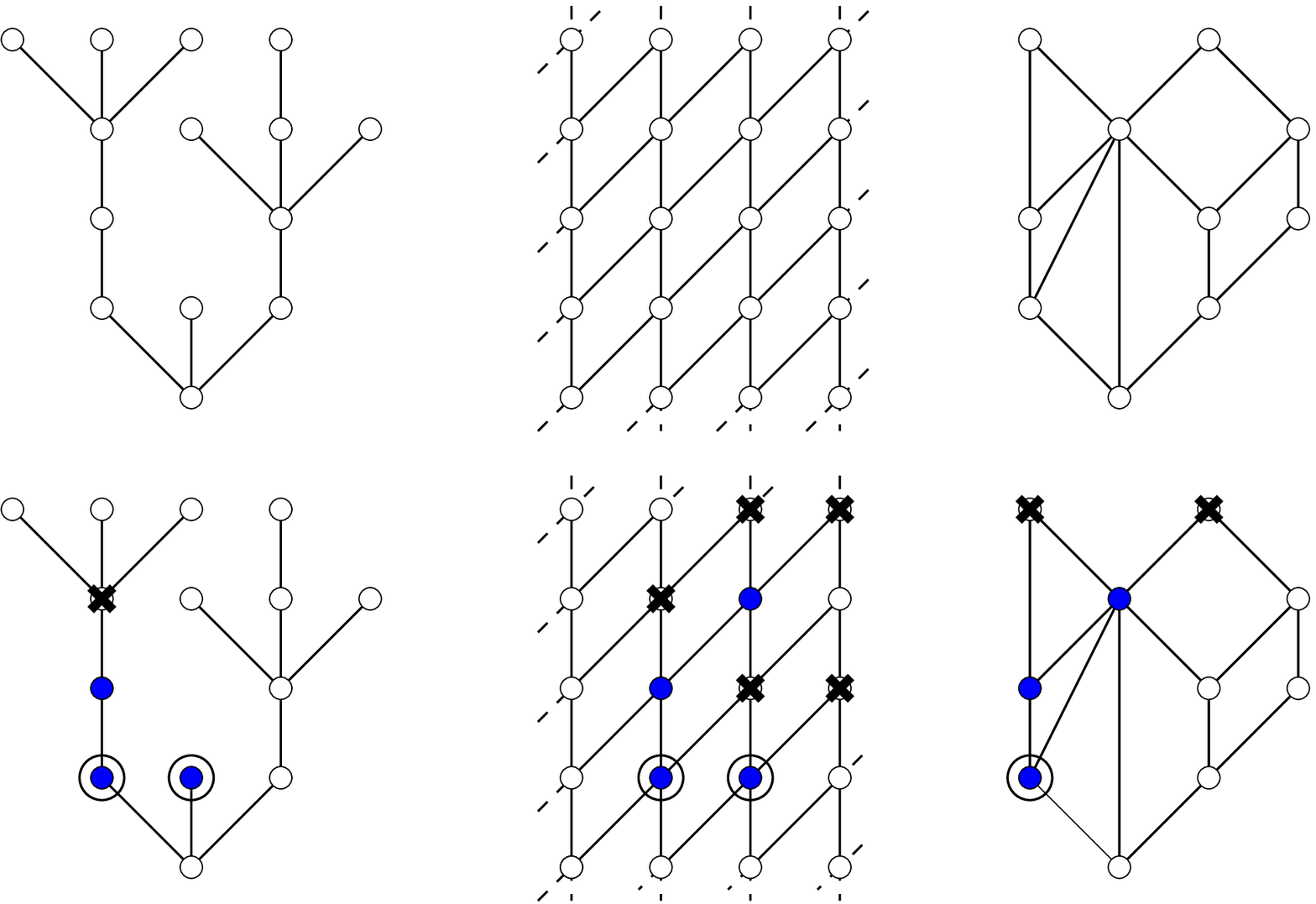}}
\captionn{\label{feffe} On the first line, three examples of agreeable graphs on which the directed edges are directed upwards. The first example is a tree, the second is the square lattice and the third a ``non-layered'' agreeable graph. On the second line are represented some DA on these graphs; filled points are the cells, the crosses are perimeter sites, and the surrounded points are the minimal sources of these DA.}
\end{figure}

We denote by ${\cal A}(S,G)$ (or more simply ${\cal A}(S)$ when no confusion on $G$ is possible), the set of finite or infinite DA on $G$ with source $S$, and by $\G^G_S$  (or $\G_S$) the generating function (GF) of finite DA counted according to the area~:
\[\G^G_S(x):=x^{|S|}+\dots\] where $\dots$ stands for a sum of monomials whose degree are at least $|S|+1$. We will need also sometimes to consider DA having their sources included (or equal) in a given set $S$. For such a DA, $S$ is called an \it over-source \rm. We denote by  $\bar{{\cal A}(S,G)}$ (or $\bar{{\cal A}(S)}$)  the set of DA having $S$ as over-source, and by $\oG^G_S$ (or $\oG_S$) its GF ($\bar{\bf G}^G_S(x)=1+|S|x+\dots$).  Finally, we set $\G^G_\varnothing=\oG^G_\varnothing=1$.

When dealing with elements $A$ of $\bar{{\cal A}(S,G)}$, the set $S\setminus A$ are also considered as (special) perimeter sites, and we set $\bar{{\cal P}_S}(A)={\cal P}(A)\cup (S\setminus A)$.

We introduce a notion of particle systems, or gas occupation, on a graph~:
\begin{defi}\label{def-gaz}Let $G=(V,E)$ be a graph. A particle system on $G$, or gas occupation on $G$, is a map $C$ from $V$ to $\{0,1\}$. The vertices $v$ such that $C(v)=1$ are said to be occupied, the others are said to be empty.
\end{defi}
In the physic literature, a hard particle model (or occupation) on a graph is a (probability model of) gas occupation $C$ with the additional constraint that the occupied sites are not neighbors in $G$. When the gas is random, we call density of $C$ in a vertex $x$, the quantity $ `P(C(x)=1)$.

\subsection{Contents}

The equivalence of enumeration of DA on lattices according to the  area and solving hard particle models has been a key point in the study of DA from its very beginning with the work of physicists, Hakim \& Nadal \cite{HV}, Nadal,  Derrida, \& Vannimenus  \cite{NDV}, Dhar \cite{DH1}, \cite{DH6} and later, thanks to some generalizations due to Bousquet-Mélou \cite{MI1} and Bousquet-Mélou \& Conway \cite{MICO}. The idea is to consider on the same lattice where are built the DA an \sl ad hoc \rm random gas models owning the same ``recursive decomposition'' as the DA up to a slight change of variables. This recursive decomposition is done using a decomposition of the lattice itself by layers (see Section \ref{SQAM}, where the historical approach used by Bousquet-Mélou \cite{MI1} is detailed).
Dhar \cite{DH5}, using statistical mechanics techniques, explained that in special cases there exists an exact equivalence between the enumeration of DA in $d$ dimensions and the computation of the free-energy of a $(d-1)$-dimensional lattice gas. A combinatorial explanation using heaps of pieces was proposed by Viennot \cite{VI1} (see also Bétréma \& Penaud \cite{BEPE} for a pedagogical and detailed exposition).\par

To solve the equation involving the GF of DA obtained by this recursive decomposition some properties of the gas model are used. The gas model is a stochastic process indexed by the lattice having in the tractable cases some nice Markovian-type properties on the layers. The rigorous arguments given in  \cite{MI1,MICO} avoid the construction of the gas model on the whole lattice as done by Dhar~: the gas model is defined on the layers of a cylinder (a lattice having some cyclical boundary conditions), and the transition allowing to pass from a layer to the following one are of Markovian type. Bousquet-Mélou \cite{MI1} finds an explicit solution of the gas process on a layer (in the square lattice case, in the triangular lattice case, and in other lattices in the joint work with Conway). Then, the computation of the density of the gas distribution is explicitly solved using that the number of configurations on such layers is finite, leading to rational GF. The GF of DA on the entire lattice (without the cyclical condition) is then obtained by a formal passage to the limit. 
\medskip

In this paper we revisit the relation between the enumeration of DA on a lattice, and more generally on any agreeable graph, and the computation of the density of a model of gas. Our construction coincides with that used by Dhar or Bousquet-Mélou in their works. 
In Section \ref{CAG}, we explain how the usual construction of random DA on a graph $G$, using a Bernoulli coloring of the vertices of $G$, allows to define in the same time a random model of gas (that we qualified to be of type 1). Here the construction is not done in ``parallel'' as in the works previously cited but on the same probability space; this provides a coupling of these objects. 
Using this coupling, we provide a general explanation of the fact that the GF of DA with source $\{x\}$ counted according to the area on any agreeable graph equals the density of the associated gas at vertex $x$, up to a simple change of variables (Theorem \ref{fond1}). This explanation is not of the same nature than in the previously cited work: the link between the density and the GF is not only formal but is explained combinatorially at the level of the DA (Section \ref{CAG}). Moreover, the construction of the gas model can be done also on (non-regular) infinite graph in a rigorous manner and no passage to the limit is needed.  

In Section \ref{Z2}, we revisit the study of DA on the square lattice; the new description of the gas model on the whole lattice allows us to provide a description of the gas model on a line (Theorem \ref{zoa}). On this line the gas model is a Markov chain which is identified. We provide then a new way to compute the GF of DA counted according to the area (Theorem \ref{zoa}). This extends to the enumeration of DA with any source on a line (Proposition \ref{geos}); till now only the case with compact sources was known. We explain also how to compute the GF of DA with sources that are not contained in a line (Remark \ref{reme}) and provide an example (Proposition \ref{geos2}).

In Section \ref{gas2}, we present an other model of gas, that we qualify to be of type 2. The density of this gas model is related to the GF of DA counted according to the area and perimeter (Theorem \ref{gfdg}). This construction explained once again at the level of object on any agreeable graph a relation used by Bousquet-Mélou \cite{MI1} in a formal way on the square lattice. Even if we haven't find any deep application to this construction, we think that it provides an interesting generic approach to the computation of the GF of DA according to the area and perimeter, and it should lead to new results in the future.

\subsubsection*{Some other references concerning DA on lattices}

One finds in the literature numerous works concerning the enumeration of DA on lattices, most of them avoids gas model considerations. We don't want to be exhaustive here (we send the reader to Bousquet-Mélou \cite{MI1}, Viennot \cite{VG,VI1} and references therein), but we would like to indicate some combinatorial works directly related to this paper. It is interesting to notice that an important part of the papers cited below are combinatorial proofs of results found before using gas techniques.\par 
First, we refer to Viennot \cite{VI1} (and  Bétréma \& Penaud \cite{BEPE}) for the algebraico-combinatorial relation between DA and heaps of pieces. This powerful point of view having some applications everywhere in the combinatorics, allows to compute the GF of DA on the triangular lattice, and by a change of variables on the square lattice. A direct combinatorial enumeration of DA on the square lattice has been done by Bétréma \& Penaud \cite{BEPE2}; they found a bijection with a family of trees :''les arbres guingois''. \par
Heap of pieces techniques have been used by Corteel, Denise \& Gouyou-Beauchamps \cite{CDG} to give a combinatorial enumeration of DA on some lattices, first counted by Bousquet-Mélou \& Conway \cite{MICO} using gas model (of type 1). 
Viennot and Gouyou-Beauchamps  \cite{GBV} provide a bijection between DA  with compact sources on the square lattice and  certain paths in the plane; they are able to enumerate these DA. Barcucci \& al. \cite{BDPP} studied DA on the square and triangular lattices with the help of the ECO method. They found some relations with permutations with some forbidden subsequences and a family of trees.

\section{Simultaneous construction of DA and gas model of type 1}
\label{CAG}
In this part, we construct on any agreeable graph $G$ a probability space on which are well-defined a model of gas -- that we qualified to be ``of type 1'' --  and a notion of random DA. The relation between the density of the gas model of type 1 and the GF of DA counted according to the area will then be explained at the level of objects on this probability space. 

\subsection{Construction of DA}
\label{COD}
Let $G=(V,E)$ be an agreeable graph. We introduce a random coloring of $V$ by the two colors $a$ and $b$.
We need to be a little bit formal here since when $G$ is infinite the existence of a probability space where such a construction is possible is not so obvious. 
We consider the probability space $\Omega=\{a,b\}^{V}$ endowed with the $\sigma$-fields ${\cal F}$ generated by the subsets having a finite cardinality (the cylinders in the usual probability terminology) and endowed by the measure product $`P_p=\l(p\delta_a+(1-p)\delta_{b}\r)^{\otimes V}$, where $\delta_a$ is the standard Dirac measure on $\{a\}$; let $C=(C_x)_{x\in V}$ be the identity map on $V$. We denote by $`E_p$ the expectation under $`P_p$.\par
In other words $C$ is a random  coloring of $V$~: under $`P_p$ the random variables $C_x$ giving the color of the vertex $x$ of $V$ are independent and take the value $a$ and $b$ with probability $p$ and $1-p$.

\begin{figure}[htbp]
\psfrag{S}{$S$}\psfrag{Sb}{$S_{\bullet}$}\psfrag{a}{$a$}\psfrag{b}{$b$}
\centerline{\includegraphics[height=3.5cm]{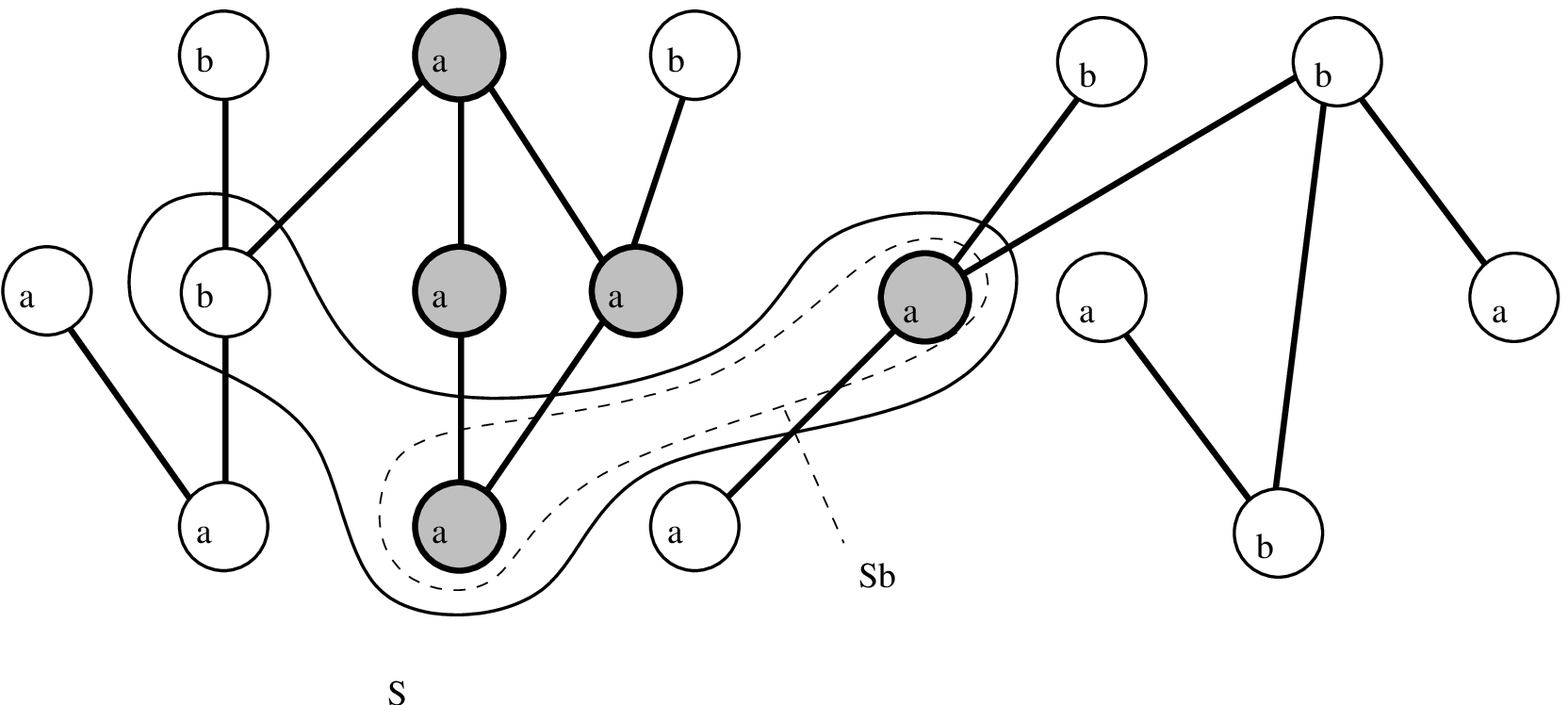}}
\captionn{\label{eeee} The DA ${\A}^{S}(`o)$ is the set of gray cells.}
\end{figure}

Let  $`o\in \Omega$ and $S$ be a subset of $V$. We denote by $S_{\bullet}(`o)=\{x, x\in S, C_x(`o)=a\}$ the subset of $S$ of vertices with color $a$. 
We denote by ${\A}^{S}(`o)$ the maximal DA for the inclusion partial order with
source $S_{\bullet}(`o)$ and whose cells are the vertices $x$ such that $C_x(`o)=a$ that can be reached from $S_{\bullet}(`o)$ by an $a$-colored path.  By construction the perimeter sites of ${\A}^{S}(`o)$ are $b$-colored (see Fig. \ref{eeee}). 

\begin{pro}\label{moch} Let $G=(V,E)$ be an agreeable graph.\par
 $(i)$ Let $S$ be a free subset of $V$ and $B\in {{\cal A}(S,G)}$ be a finite DA with source $S$. We have
\[`P_p({\A}^{S}=B)=p^{|B|}(1-p)^{|{\cal P}(B)|}.\]
$(ii)$ Let $B\in \bar{{\cal A}(S,G)}$ be a finite DA with over-source the free set $S$. We have
\[`P_p({\A}^{S}=B)=p^{|B|}(1-p)^{|\bar{{\cal P}_S}(B)|}.\]
\end{pro}
\Proof $(i)$ First ${\A}^{S}(`o)=B$ if and only if $C_x=a$ for all $x\in B$ and if $C_x=b$ for all $x\in {\cal P}(B)$. Then  $\{{\A}^{S}=B\}$ is ${\cal F}$--measurable since it depends only on a finite number of cells, and the conclusion follows. For $(ii)$, we use moreover that $\bar{{\cal P}_S}(A):={\cal P}(A)\cup ( S\setminus A)$.~$\Box$

\subsection{Directed animals and percolation}
When $G$ is an infinite graph and $|S|\geq 1$, under $`P_p$ the random DA ${\A}^{S}$ may be infinite with positive probability. 
The probability to have an infinite DA with source $S$  is also that of the directed sites percolation starting from $S$ where the cells of the percolation cluster are the vertices with color $a$ reachable from $S$ by an $a$-colored directed path. \par
 Denote by $p_{crit}^S$ the threshold for the existence of an infinite DA with positive probability~:
\[p_{crit}^S=\sup\{\,p,~ `P_p(|{\A}^S| <+\infty)=1\}.\]
Most of the results of the present paper are valid only when $p<p_{crit}^S$. The threshold $p_{crit}^S$ is in general difficult to compute, but here is a simple sufficient condition on $G$ for which $p_{crit}^S>0$.
\begin{pro}\label{yopop}
Let $G$ be a agreeable graph such that the maximum number of children of its vertices is bounded by $K$. Then for any finite subset $S$ of $G$, $p_{crit}^S\geq 1/K$.
\end{pro}
A proof of that result is given in the Appendix. 
Also in the Appendix, Comment \ref{soupe} provides a graph in which $p_{crit}=0$. \medskip

We recall two results giving some insight on the percolation probabilities (and easy to prove).\\
$\bullet$ Let $S_1$ and $S_2$ be two subsets of $V$. We have $p_{crit}^{S_1\cup S_2}=\min \left\{\,p_{crit}^{S_1}\, ,\, p_{crit}^{S_2}\right\}.$\\
$\bullet$ Let ${\bf p}_{crit}=\inf_{v\in V} p_{crit}^{\{v\}}$. 
For any $p \in[0, {\bf p}_{crit})$, under $`P_p$,  almost surely (a.s.) all DA in $G$ having a finite source are finite. 
Indeed a countable intersection of a.s. events is a.s..

\subsection{Construction of the gas model of type 1}
\label{cgm}
The construction of the gas model of type 1, Proposition \ref{tru} and the Nim game construction presented in this section are generalizations and formalization of the work of the first author \cite[section 1.4]{LB}. \medskip

Let us build a gas model $X$ on an agreeable graph $G=(V,E)$ (see Definition \ref{def-gaz}). This construction takes place on the probability space $\Omega$ introduced in Section \ref{COD}, and $X$ is defined thanks to the random coloring $C$. \par
For any $x\in V$ and $`o\in \Omega$, denote by 
\ben
\label{yop2}
X_x(`o)&:=&
\left\{
\begin{array}{ll}
0 & \textrm{ if } C_x(`o)=b\\
\dis\prod_{c : \textrm{ children of }x} (1-X_{c}(`o))& \textrm{ if }C_x(`o)=a
\end{array}\right.\\
&=& \1_{C_x(`o)=a} \dis\prod_{c : \textrm{ children of }x} (1-X_{c}(`o)).
\een
If $x$ has no children the product in \eref{yop2} is empty, and as usual, we set its value to 1.\par 
For any $`o$ and any $x\in V$, $X_x(`o)$ is to be interpreted as the gas occupation in the vertex $x$. \medskip

\begin{figure}[htbp]
\psfrag{a}{$a$}\psfrag{b}{$b$}
\centerline{\includegraphics[height=7 cm]{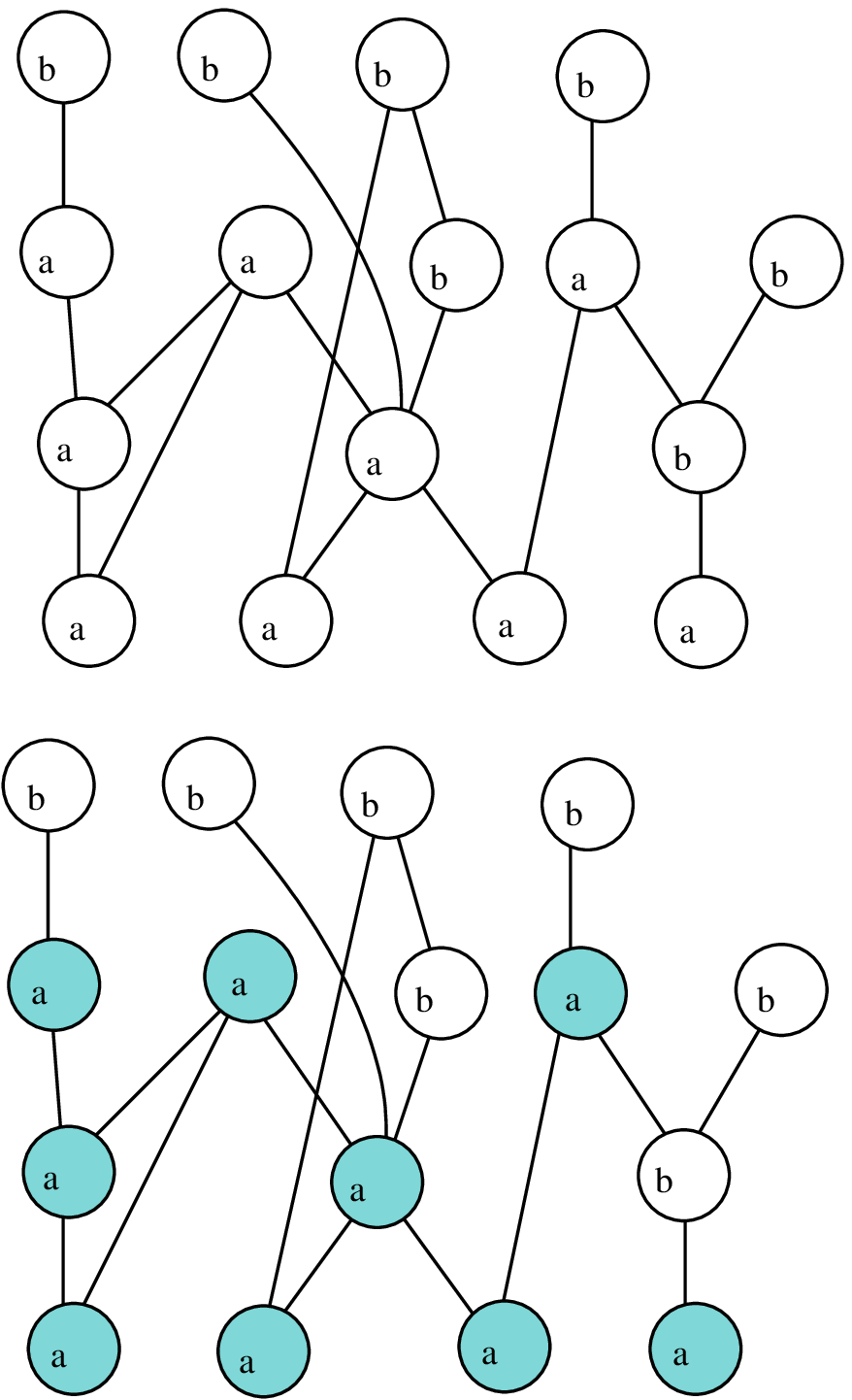}~~~\psfrag{b}{0}\psfrag{1}{1}\psfrag{a}{$\star$}\psfrag{0}{$0$}\includegraphics[height=7 cm]{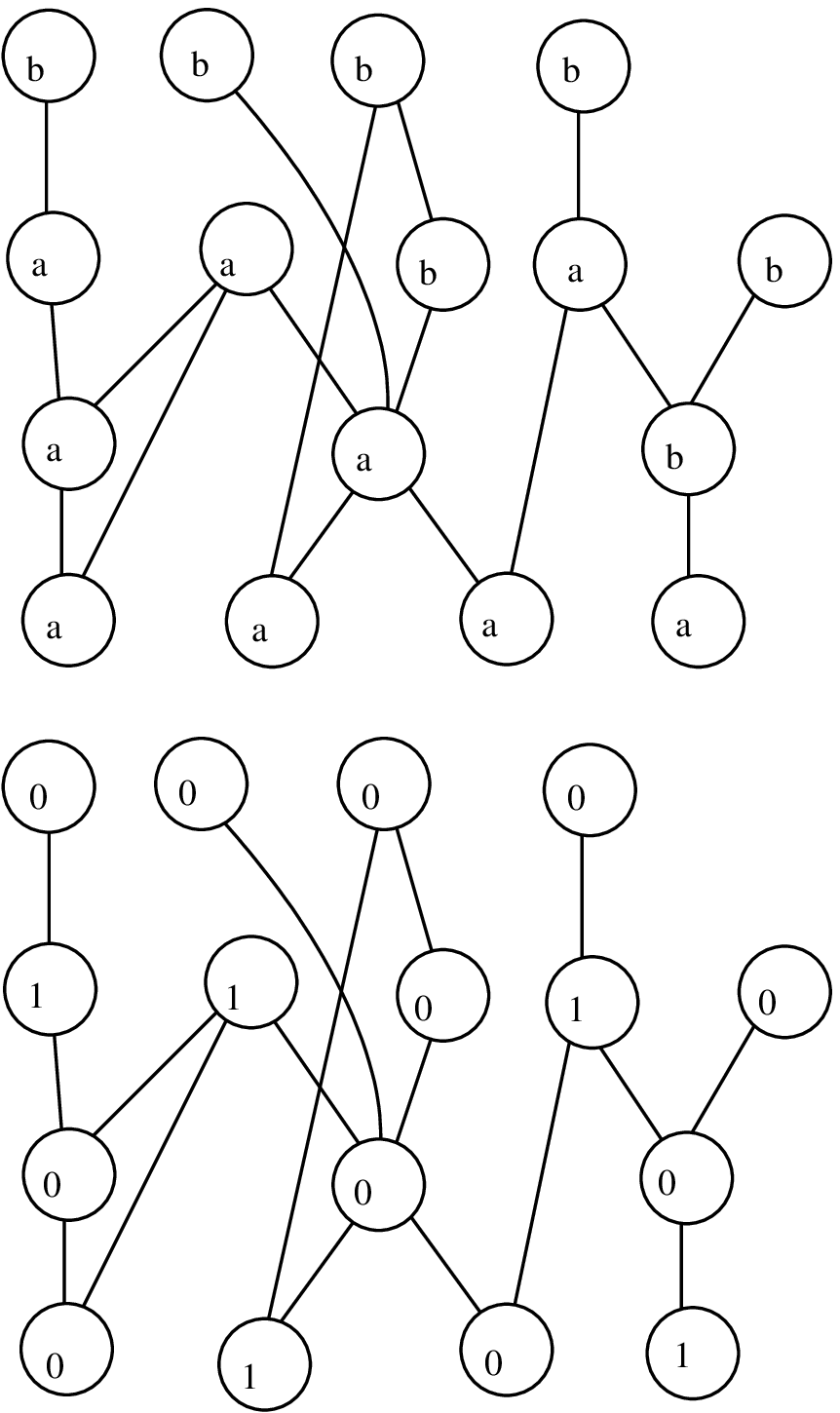}}
\captionn{\label{tir} On the first column on top, a random coloring. Below, the family of DA derived from it. On the second column on top, the beginning of the computation of the gas occupation. ``$\star$'' stands for the places where the calculus $X_x\leftarrow \prod_{c : \textrm{ children of }x} (1-X_{c})$ must be done. Below, the gas occupation has been computed.}
\end{figure}

We have to investigate when the recursive definition giving $X_x(`o)$ is correct, that is when it allows to indeed compute a value $X_x(`o)$ (see Fig. \ref{tir}). \\
$\bullet$ When $C_x(`o)=b$ then $X_x(`o)=0$ (there is no problem to define $X_x(`o)$).\\
$\bullet$ When $C_x(`o)=a$, to compute $X_x(`o)$, knowing all the values $X_y(`o)$ for $y$ child of $x$ is sufficient; their values are given by the same rule.
By successive iterations a ``DA of calculus'' growths. This DA is exactly ${\A}^{\{x\}}(`o)$.  Indeed, on each cell of ${\A}^{\{x\}}(`o)$ -- whose color are $a$ by construction -- the following computation is done 
\[X_x(`o) \leftarrow \dis\prod_{c : \textrm{ children of }x} (1-X_{c}(`o));\]  since $C_y(`o)=b$ for any perimeter sites of ${\A}^{\{x\}}(`o)$,  $X_y(`o)=0$ on ${\cal P}({\A}^{\{x\}})(`o))$.
Then, one sees that $X_x(`o)$ is well defined if ${\A}^{\{x\}}(`o)$ is finite because this recursive computation of $X_x(`o)$ ends.   \par
Moreover, in this case, if  ${\A}^{\{x\}}(`o)=A$ is finite  the value $X_x(`o)$ is a deterministic function of $A$, that we denote by $\chi_x(A)$, since the only data in this calculus is the geometry of $A$ (the map $\chi_x$ is defined only on finite DA with source $x$). 
The maps $(\chi_x)_{x\in V}$ satisfies then for simple reasons the  following decomposition.
Let $v$ be a vertex in $G$, $A^{\{v\}}$ a finite DA with source $v$, and denote by $v_1,\dots, v_d$ the children of $v$ in $G$, and  ${A}^{\{v_1\}},\dots, {A}^{\{v_d\}}$ be the maximal DA included in $A^{\{v\}}$ with over-source $\{v_1\},\dots,\{v_d\}$ respectively.
Then 
\begin{equation}\label{mort}
\chi_v\l({A}^{\{v\}}\r)=\1_{|{A}^{\{v\}}|>0}\prod_{i=1}^d\l(1-\chi_{v_i}({A}^{(v_i)})\r).
\end{equation}
In the same vein, assume that a finite free subset $S$ of $G$ is given. The vector $(X_x(`o))_{x\in S}$ giving the gas occupation on $S$ is also a deterministic function of the DA ${\A}^{S}(`o)$.\par
\begin{remi}In the case where $C_x(`o)=a$ but $|{\A}^{\{x\}}(`o)|$ is infinite, the computation of $X_x(`o)$ may also ends within a finite number of steps (using the fact that the product $\prod_{c : \textrm{ children of }x} (1-X_{c}(`o))$ is known to be 0 when one of its terms is null, which can be the case on a finite sub animal of ${\A}^{S}(`o)$). But, in this case the value of $X_x$ returned by this procedure is inadequate and must be not considered in the analysis since the a.s. finitness of ${\A}^{S}$ is crucial in the proofs. In the following we will restrict ourselves to $p<p_{crit}$ in order to avoid infinite DA ${\cal A}^S(`o)$ with probability 1.
\end{remi}

By the previous consideration we may conclude by the following proposition.
\begin{pro}\label{tru}Let $G=(V,E)$ be an agreeable graph, 
$x\in V$ and $p\in[0,p_{crit}^{\{x\}})$. Under $`P_p$ the random variable $X_x$ is a.s. well defined by \eref{yop2}, and
\[`E_p(X_x)=`P_p(X_x=1)=`E_p(\chi_x({\A}^{\{x\}}))=\sum_{A\in {\cal A}(\{x\},G), |A|<+\infty, \chi_x(A)=1} p^{|A|}(1-p)^{|{\cal P}(A)|}.\]
where ${\cal A}(\{x\},G)$ has been defined in Section~\ref{gc}.
\end{pro}

\proof For $p\in[0,p_{crit}^{\{x\}})$, the DA ${\A}^{\{x\}}$ is a.s. finite, and then a.s. $\chi_x({\A}^{\{x\}})=X_x$. The three equalities are straightforward. One also checks that $X_x$ is indeed a random variable, that is ${\cal F}$-measurable. For $i\in\{0,1\}$, 
\[X_x^{-1}(\{i\})=\bigcup_{n\geq 0}\{A, |A|=n, \textrm{ DA with source }\{x\}, \chi_x(A^{\{x\}})=i\}.\] This sets are measurable since they depend only on a finite number of cells. The conclusion follows the fact that a countable union of measurable sets is measurable. 
~$\Box$ 

\begin{remi}
In the case where ${\A}^{\{x\}}(`o)=A$ is finite, it may be convenient for the reader to see the computation of $X_x(`o)$ as the result of a Nim type game with two players, player(0) and player(1) on ${A}$ according to the following rules~: \\
-- the first player that can not play is the looser.\\
-- at time 0,  player(0) places a token on $x$ (if  $A$ is empty, $player(0)$ is the looser). \\
-- then, the players move in turn the token upward in $A$. At time $i$, $player(i\mod 2)$ takes the token from where it is, say in $v$, and can move  it in $u$, if $u$ is in $A$ and if $(v,u)$ is a directed edge of $G$.  \par
The fact is that $\chi_x(A)=1$ iff player(0) has a strategy to win against all defense of player(1). We let the reader proves this property as an exercise. Indication: proceed to the computation of $\chi_x(A)$ by the last moves of the players.
\end{remi}

One of the main aim of this paper is to provide an explication at the level of objects of the following theorem 
\begin{theo}\label{fond1}
Let $G=(V,E)$ be an agreeable graph, $x$ be  a vertex of $V$, $R_{x}$ be the radius of convergence of $\G^G_{\{x\}}$ and $p\in[0,R_x)$. 
We have
\begin{equation}\label{fond2}
`E_p(X_x)= -\G^G_{\{x\}}(-p).
\end{equation}
\end{theo}
We will see in Lemma \ref{lemRx} that in any graph $R_{x}\leq p_{crit}^{\{x\}}$.\medskip

This relation between the density of a gas model and the GF of DA in the case of lattice graph is first discovered by Dhar \cite{DH5} in the square lattice case; it is then made rigorous and generalized by Bousquet-Mélou \cite{MI1,MICO}. In each case, a formula similar to \eref{fond2} is obtained, but there the equality is only formal~: $`E_p(X_v)$ and $-\G^G_{\{v\}}(-p)$ are shown to be formal series satisfying the same recursive decomposition. We want also to point out that in \cite{MI1,MICO}, the gas model is studied on a cylinder and some arguments using the finiteness of the number of states, and the convergence of some Markov chains are used (see Section \ref{SQAM}). Here these steps are not needed since the construction of the gas model is done ``on the right graph at once'' and the construction provides also the gas under ``the stationary regime'' if the graph is a lattice. We want also to stress on the fact that Theorem \ref{fond1} holds on any agreeable graph and not only on lattices.
\par

In order to express our relation in the level of objects, we will rewrite the right hand side of \eref{fond2} under the form of an expectation, in order to make more apparent that $`E_p(X_x)$ and $-\G^G_{\{x\}}(-p)$ are both sum on weighted DA (the weight of a DA $A$ being simply $p^A(1-p)^{|{\cal P}(A)|}$ its probability on $\Omega$). But first, we need to introduce the notion of sub animal.

\begin{defi} Let $A$ and $A'$ be two DA. We say that $A$ is a sub animal (sub-DA) of $A'$ (we write $A\prec A'$), if $A$ is included in $A'$ and if ${\cal S}(A)$ is included in ${\cal S}(A')$. 
\end{defi}

For any DA $A$ with minimal source ${\cal S}(A)=S$, and $p\in(0,p_{crit}^S)$, we have
\[p^{|A|}=\sum_{A'~:~ A\prec A',{\cal S}(A')=S}`P_p({\A}^S=A')=\sum_{A'~:~ A\prec A',{\cal S}(A')=S}p^{|A'|}(1-p)^{|{\cal P}(A')|}.\]
Indeed, $p^{|A|}$ is the probability of the event $\{`o,X_x(`o)=a \textrm{ for any }x\in A\}$, and this latter equals $\{`o, A \textrm{ is a sub-DA of }{A}^S(`o) \}$. 
We have

\begin{pro}\label{etap}  Let $G=(V,E)$ be an agreeable graph.
For any $x\in V$ and any $p\in[0,R_x)$,
\[-\G^G_{\{x\}}(-p)=`E_p\l(D_x({\A}^{\{x\}})\r)\]
where for any $A$ with source $\{x\}$,
\begin{equation} \label{diff}
D_x(A)= \sum_{B~:~ B \prec A, {\cal S}(B)=\{x\}} (-1)^{|B|+1}
\end{equation} is the difference between the number of sub-DA of $A$ having an odd number of cells and those having an (non zero) even number of cells. 
\end{pro}
Before proving this proposition, we establish the following Lemma
\begin{lem}\label{lemRx} For any agreeable graph $G=(V,E)$ and $x$ in $V$,
\begin{equation}\label{ineq}
R_{x}\leq p_{crit}^{\{x\}}.
\end{equation}
\end{lem}
\proof Since for any $A$ with source $\{x\}$, $p^{|A|}=`P(A \textrm{ is a sub-DA of } {\A}^{\{x\}})$ (where ${\A}^{\{x\}}$ may be infinite),
\begin{equation}\label{rc}
\sum_{A} p^{|A|}=\sum_A  `P(A \prec {\A}^{\{x\}})=\sum_A `E_p(1_{A \prec {\A}^{\{x\}}})= `E_p(\textrm{ Number of sub-DA of }{\A}^{\{x\}}). 
\end{equation}
If $p\in[0,R_{x})$ then these quantities are finite. This implies that $`P_p$ a.s. the number of sub-DA of ${\A}^{\{x\}}$ is finite, which in turn implies that ${\A}^{\{x\}}$ is $`P_p$ a.s. finite, which finally yields $p\leq p_{crit}^{\{x\}}$. ~$\Box$ \medskip

\noindent \bf Proof of Proposition \ref{etap}\rm. It relies on a permutation of sum symbols. 
For any $p\in[0,R_{x})$, 
\be
-\G^G_{\{x\}}(-p)&=&\sum_{A, {\cal S}(A)=\{x\}} p^{|A|}(-1)^{|A|+1}.
\ee
Since by Lemma \ref{lemRx} we have $p<p_{crit}^{\{x\}}$, this is equal to 
\begin{equation}\label{zrio}
\sum_{A~:~ {\cal S}(A)=\{x\}}\l(\sum_{A'~:~A\prec A', {\cal S}(A')=\{x\}}p^{|A'|}(1-p)^{|{\cal P}(A')|}(-1)^{|A|+1}\r). 
\end{equation}
This sum converges absolutely when $p<R_{x}$ since with an absolute value it is the mean of the number of Sub DA of ${\A}^{\{x\}}$ (which is finite when $p<R_{x}$, see the proof of Lemma~\ref{lemRx}). By the Fubini's theorem the double sum in \eref{zrio}  is
\be
\sum_{A'~:~{\cal S}(A')=\{x\}}p^{|A'|}(1-p)^{|{\cal P}(A')|} \l(D_x(A')\r)
=`E_p\l(D_x({\A}^{\{x\}})\r).~\Box
\ee

\begin{remi} 
$\bullet$ In general the inequality \eref{ineq} is strict, since 
\[R_x=\sup\{p, `E_p(\textrm{Number of sub-DA of }{\A}^{\{x\}}) <+\infty\} 
 \textrm{~~when~~} p_{crit}=\sup\{p, `P_p(|{\A}^{\{x\}}|<+\infty)=1\}.\] In the case of the square lattice $R_x=1/3$ when $p_{crit}$ is larger than $1/2$ (by Proposition \ref{yopop}) and is expected to be around 0.6.\\
$\bullet$ Proposition \ref{etap} is a relation between two quantities that are not defined in general for the same values of $p$. Viewed as formal series in $p$, they are equal. But viewed as real functions of $p$, this is not so simple, and this is rather interesting. This is due to the classical well known fact that a sum which is not absolutely converging is well defined only if an order of summation of its terms is provided. When $p>R_x$ the fact that $-\G^G_{\{x\}}(-p)$ is a non convergent sum comes from the fact that the order of summation is given by the sizes of the DA. The second sum $\sum_{A} p^{|A|}(1-p)^{|{\cal P}(A)|} \l(D_x(A)\r)$ contains the same terms as $-\G^G_{\{x\}}(-p)$ as one may guess by the permutation of the sum symbols authorized for some $p$ by the Fubini's theorem. But these terms appear in some sense by groups, and it is a fact that this sum converges more easily  than $-\G^G_{\{x\}}(-p)$ (as said in the first point).
\end{remi}

We saw that $-{\G^G_{\{x\}}}(-p)=`E_p\l(D_x({\A}^{\{x\}})\r)$ for $p<R_{x}\leq p_{crit}^{\{x\}}$ and the density of the gas $`E_p(X_x)$ is equal to $`E_p({\chi_x({\A}^{\{x\}})})$, for $p<p_{crit}^{\{x\}}$. Here is the explication of Theorem \ref{fond1} at the level of object~:
\begin{theo} \label{zaea}For any finite DA $A$ with source $v$, we have
$D_v(A)= \chi_v(A)$. \end{theo}
This allows to deduce that Theorem \ref{fond1} holds true, since a.s. when $p<p_{crit}^{\{v\}}$, the random DA ${\A}^{\{v\}}$ is a.s. finite, and then $X_v({\A}^{\{v\}})=\chi_v({\A}^{\{v\}})=D_v({\A}^{\{v\}})$ a.s., and then these variables have the same expectation (notice that this also implies that $D_v(A)$ takes its values in $\{0,1\}$, which is not necessarily obvious). \par
 
In order to prove Theorem \ref{zaea}, we will show that $D_v$ owns the same recursive decomposition as $\chi_v$ given in \eref{mort}; this is done via the introduction of a notion of embedding of trees in DA. An heuristic is given in Section \ref{zoza}.

\subsection{Embedded trees}

Let $G=(V,E)$ be an agreeable graph. 
The set $V$ being at most countable we assume from now on that an order denoted by $\sous{<}{~V}$ is given on $V$. This order induces an order among the children of a given vertex in $G$.
Since we are to ``canonically'' embed some ordered trees in $G$ we need also to define a suitable ordering of the nodes of those trees; this is inspired from the Neveu's definition of trees. \par 

\begin{figure}[htbp]
\psfrag{(2,6)}{$2,6$}\psfrag{B}{$b$}\psfrag{C}{$c$}
\psfrag{(1,2)}{$1,2$}
\psfrag{(2,5)}{$2,5$}
\psfrag{(7,5)}{$7,5$}
\psfrag{(1,1)}{$1,1$}
\psfrag{1}{$1$}\psfrag{7}{$7$}
\psfrag{2}{$2$}\psfrag{e}{$\varnothing$}
\psfrag{3}{$3$}\psfrag{4}{$4$}\psfrag{(3,1)}{$3,1$}
\centerline{\includegraphics[height=2cm]{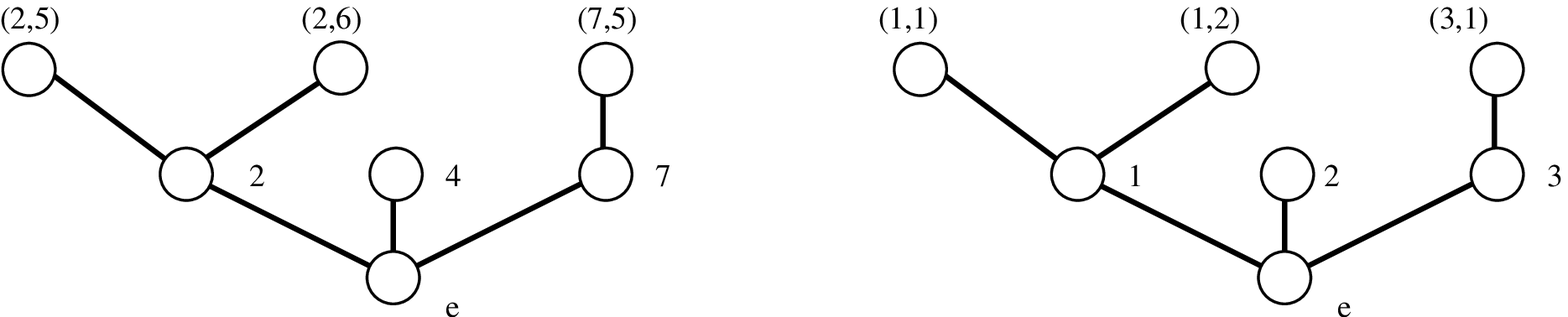}}
\captionn{\label{nev} The trees admit some usual representation as embedded figure in the plane, mapping the order between brothers into the usual order between the abscissas. In the usual Neveu's convention, there is an additional axiom (if $v=u_1\dots u_kj\in \bt$ for some $j\geq 2$ then $v=u_1\dots u_k(j-1)\in \bt$) which ensures that the ``branching structure of the tree'' characterizes the tree (for Neveu, only the second drawing represents a tree). This is not the case, here. }
\end{figure}

Let $`N=\{1,2,3,\dots\}$ be the set of non-negative integers and ${\cal W}=\{\varnothing\}\cup \dis\bigcup_{i=1}^{+\infty}`N^i$ the set of finite words on the alphabet $`N$, where $\varnothing$ denotes the empty word. We define the concatenation product of two words $u=u_1\dots u_k$ and $v=v_1\dots v_l$ of ${\cal W}$ by  $uv:=u_1\dots u_kv_1\dots v_l$; the empty word $\varnothing$ is the neutral element for this operation~: $\varnothing u=u\varnothing$ for any $u\in{\cal W}$.
\begin{defi}
We denote by tree a subset $\bt$ of ${\cal W}$ such that 
$\varnothing\in\bt$ and if $u=u_{1}\dots u_k\in \bt$ for some $k\geq 1$ then $u_{1}\dots u_{k-1}\in \bt$. In other word if a word $u$ is in $\bt$, its prefixes are also in $\bt$ (see Fig. \ref{nev}).
\end{defi}
The set of trees is denoted by ${\cal T}$.\par 
In the combinatorial literature, what is called tree here is sometimes viewed as a depth first traversal encoding of trees.\par
The vocabulary attached to trees are as usual: $\varnothing$ is called the root, the strict  prefixes of $u\in\bt$ are called ancestors of $u$, and $u$ is a descendant of its ancestors. The elements of $\bt$ are called nodes, and for $u\in\bt$, $|u|$ stands for the number of letters in $u$ (by convention $|\varnothing|=0$) and is called the depth of $u$. If $u=u_1\dots u_ki$ then $u_1\dots u_k$ is the father of $u$, and if $v=u_1\dots u_kj\in \bt$ for $i\neq j$, we say that $u$ and $v$ are brothers. For any $u=u_1\dots u_k\in\bt$, we denote by ${\cal C}_u(\bt)=\{ui,i\in`N\}\cap \bt$ the set of children of $u$ in $\bt$. The size of $\bt$ denoted by $|\bt|$ is the cardinality of $\bt$.

\begin{figure}[htbp]
\psfrag{v}{$\varnothing$}
\psfrag{1}{1}
\psfrag{w}{$v$}
\psfrag{3}{3}
\psfrag{31}{3,1}
\psfrag{33}{3,3}
\psfrag{4}{4}
\psfrag{6}{6}
\psfrag{11}{1,1}
\psfrag{14}{1,4}
\psfrag{143}{1,4,3}
~\psfrag{43}{4,3}
\psfrag{62}{6,2}
\psfrag{431}{4,3,1}
\psfrag{433}{4,3,3}
\centerline{\includegraphics[height=3.3 cm]{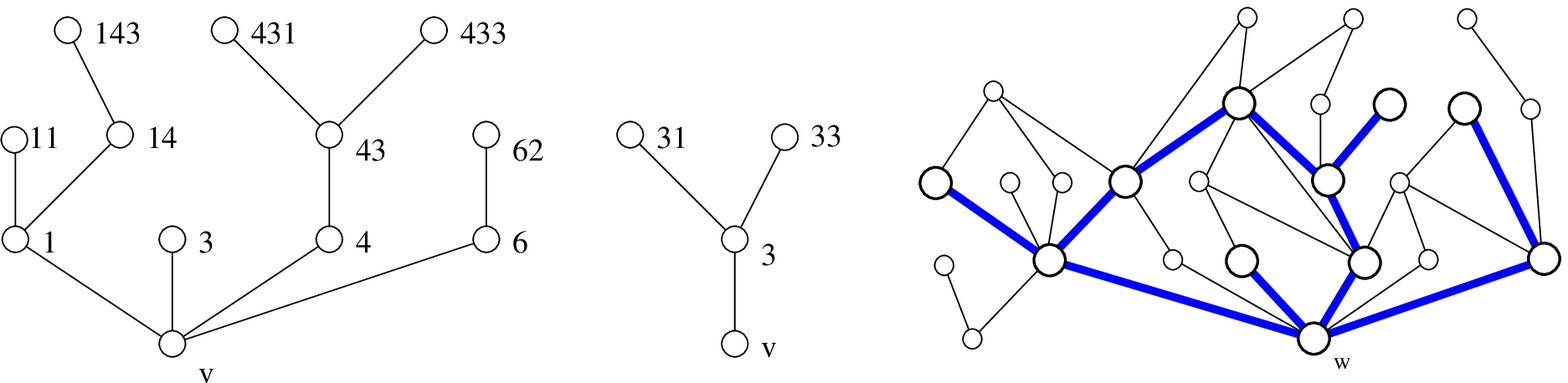}}
\captionn{\label{eaih} A representation of the tree $\bt=\{\varnothing,1,3,4,6,(1,1),(1,4),(4,3),(6,2),(1,4,3),(4,3,1),(4,3,3)\}$. 
The second tree is $\bt^4$. We have ${\cal C}_\varnothing(\bt)=\{1,3,4,6\}$. On the third figure is drawn the embedding $\pi_v(\bt)$ on an agreeable graph, where the edges are directed upward.}
\end{figure}

Let $\bt$ be a tree and $i\in{\cal C}_\varnothing(\bt)$ a child of the root in $\bt$. For any $i\in\mathbb{N}$, we denote by $\bt^i$ the set~:
\[ \bt^i:=\{ v, u:=iv \in \bt\}.\]
It is easy to see that $\bt^i$ is a tree, and that it is obtained from $\bt$ by deleting all the nodes that are not descendant of $i$, and ``rerooted'' in $i$~: $i$ and the descendants of $i$ are kept but they lose they first letter in order to form a tree (see Fig. \ref{eaih}). The tree $\bt^i$ is  called \it fringe subtree \rm in the literature.
 
\begin{defi}Let $G=(V,E)$ be an agreeable graph, $v\in V$ and $\bt$ a tree. Let $v_1,\dots,v_d$ be the $d=d_v$ children of $v$ in $G$ sorted according to $\sous{<}{~V}$. We say that $\bt$ is embeddable in $G$ at $v$ if $\bt=\{\varnothing\}$ or if for any $i\in {\cal C}_\varnothing(\bt)$, $\bt^i$ is embeddable at $v_i$. If $\bt$ is embeddable in $G$ at $v$, we denote by $\pi_v(\bt)$ its embedding in $G$ at $v$ defined by~:
\[\pi_v(\bt)=\{v\}\cup\bigcup_{i\in{\cal C}_\varnothing(\bt)}\pi_{v_i}(\bt^i).\]
An illustration of this embedding is given in Fig. \ref{eaih}.
\end{defi}
In other words, $\pi_v(\bt)$ is the subset of $G$ obtained by drawing $\bt$ on $G$ according to the following rules:~ first draw the root of $\bt$ on $v$ (in other word $\pi_v(\{\varnothing\})=v)$. Each branch of $\bt$ is a succession of nodes $a_k=u_1\dots u_k$ where $a_k$ is a prefix of $a_{k+1}$ (and $a_0=\varnothing$). Then draw $a_{k+1}$ in $G$ in such a way that $a_{k+1}$ is the $u_{k+1}$ child of $a_k$ (the edge $(\pi_v(a_k),\pi_v(a_{k+1}))$ of $G$ is the $u_{k+1}$th edge starting from $\pi_v(a_k)$) (see Fig. \ref{eaih}). Hence each branch of $\bt$ is finally embedded in a simple path of $G$ issued from $v$.\medskip

For any tree $\bt$ embeddable in $G$ at $v$, $\pi_v(\bt)$ is a DA on $G$ with source $\{v\}$.
For any finite DA $A^{\{v\}}$ with source $\{v\}$ in $G$, we set 
\begin{equation}\label{Delta}
\Delta_v(A^{\{v\}}):=\sum_{\bt~:~\pi_v(\bt)\prec A^{\{v\}}} (-1)^{|\bt|+1},
\end{equation}
where the sum is taken on the (necessarily finite) set of trees embeddable in $G$ at $v$. %
Our proof of Theorem~\ref{zaea} consists of proving the two following equalities valid for any finite DA $A$ with source $\{v\}$:
\begin{equation}\label{dix}
 D_{v}(A) =  \Delta_v(A) \mbox{ and }  \Delta_v(A) = \chi_v(A). 
\end{equation}

The key lemma in the proof of the first equality is 

\begin{figure}[htbp]
\psfrag{A}{$a$}\psfrag{B}{$b$}\psfrag{C}{$c$}
\psfrag{v}{$\varnothing$}\psfrag{vv}{$v$}
\psfrag{w}{$w$}
\centerline{\includegraphics[height=4cm]{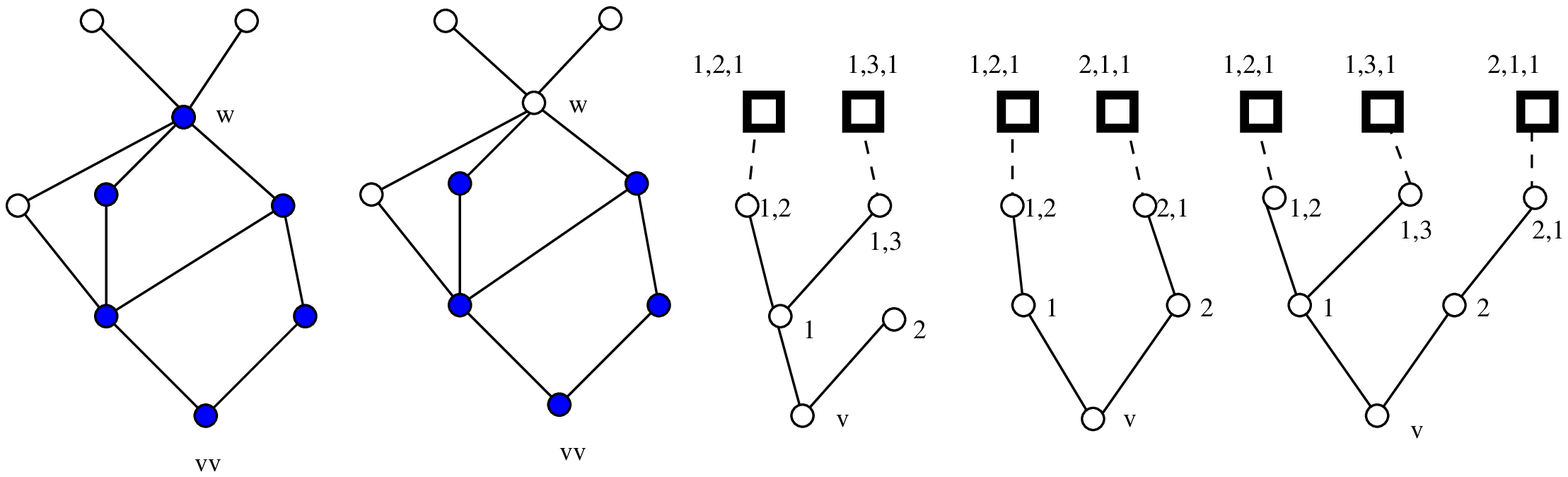}}
\captionn{\label{inv} On the first picture, a DA $A^{\{v\}}$ (whose cells are colored). On the second one, a DA $A'$ obtained from $A^{\{v\}}$ by the suppression of a cell having no child. The three trees with straight lines represent trees embeddable in $A'$ at $v$. The trees with in addition some of the doted lines are the trees embeddable in $A^{\{v\}}$ at $v$. The square nodes are the nodes whose embedding are $w$}
\end{figure}

\begin{lem} \label{a} For any finite DA $A^{\{v\}}$ with source $\{ v\}$, we have 
\[ \sum_{\bt~:~ \pi_v(\bt)= A^{\{v\}}} (-1)^{|\bt|} = (-1)^{|A^{\{v\}}|}\]
\end{lem} 
{\bf Proof :} This is true if $A^{\{v\}}=\{v\}$.
Assume that this is true for all DA having $k$ cells or less than $k$ cells. Take a DA $A^{\{v\}}$ with $k+1$ cells. Among those cells, there exists $w$ such that $w$ has no child  in $A^{\{v\}}$. %
Consider the DA $A'=A^{\{v\}}\setminus\{w\}$ and compare the set of trees $\xi$ whose embedding is $A^{\{v\}}$ with $\xi'$ the set of trees whose embedding is $A'$. The trees in $\xi$ are obtained from those of $\xi'$ by addition of leaves which embedding is $w$. Hence by the natural  projection from $\xi$ into $\xi$', (see Fig. \ref{inv})
\begin{equation}\label{rose}
\xi=\bigcup_{t\in \xi'} K_t
\end{equation}
where $K_t$ is the subset of $\xi$ containing the trees obtained from $t$ by the addition of some leaves. The $K_t$'s form a partition of $\xi$. 
One has, for any $t\in \xi'$,
\[\sum_{\tau \in K_t} (-1)^{|\tau|} =\sum_{\tau \in K_t} (-1)^{|t|+|\tau\setminus t|}.\]
 But given $t$, there is a maximal tree in $K_t$ (for the inclusion) obtained from $t$ by the addition of a set $R_t$ of leaves. The others trees of $K_t$ are formed by $t$ together with a non-empty part of $R_t$. Hence
\[\sum_{\tau \in K_t} (-1)^{|\tau|} =\sum_{B, B\subset R_t, |B|\geq 0}  (-1)^{|t|+|B|}=(-1)^{|t|}((1-1)^{|R_t|}-1)=(-1)^{|t|+1}.\]
By \eref{rose}, $\sum_{\tau \in \xi} (-1)^{|\tau|}=-\sum_{t\in \xi'} (-1)^{|t|}$ and we conclude by recurrence on the number of cells of $A$.
~$\Box$
\medskip

Now, we conclude the proof of  Theorem~\ref{zaea}.\\
{\bf Proof of  Theorem~\ref{zaea}. }
On one hand, by Lemma~\ref{a}, $D_v(A)$ (defined in \eref{diff}) satisfies 
\be D_v(A)=  \sum_{B : B \prec A, {\cal S}(B) = \{ v \}}\left( \sum_{\bt~:~ \pi_v(\bt)= B} (-1)^{|\bt|+1}\right) = \sum_{\bt~:~ \pi_v(\bt) \prec A} (-1)^{|\bt|+1} = \Delta_v(A).
\ee
We now establish that $\Delta_v$ owns the same recursive decomposition (\ref{mort}) as $\chi_v$. If $|A|=1$ then $\Delta_v(A) = 1 = \chi_v(A)$. 
Assume that the source $v$ of $A$ has $d$ children $v_1,\dots,v_d$. For any $i\in\{1,\dots,d\}$, we have
\[\Delta(A^{v_i}):=-\sum_{\bt~:~ \pi_{v_i}(\bt)\prec A^{{v_i}}} (-1)^{|\bt|},\]
where $A^{v_i}$ is the maximal DA included in $A$ with source $v_i$.
Decomposing the trees $\bt$ embeddable in $A$ according to their embeddings in the $A^{v_i}$, we get
\be
\Delta_v(A^{v}) & = & 1_{|A^v|>0}\l(\sum_{B\subset \{v_1,\dots,v_d\}} \prod_{i\in B}(-1)^{|\bt^i|}\r)\\
           & = & 1_{|A^v|>0}\l(\sum_{B\subset \{v_1,\dots,v_d\}} (-1)^{|B|}\prod_{i\in B}\Delta_{v^i}(A^{{v_i}})\r)
=1_{|A^v|>0} \prod_{i=1}^d( 1-\Delta_{v^i}(A^{{v_i}})),
\ee
Hence $\Delta_v$ owns the same recursive decomposition (\ref{mort}) as $\chi_v$, and finally $\Delta_v$ and $\chi_v$ coincide on the set of finite DA with source $v$.
~$\Box$ 

\subsection{Why is the embedding of trees in DA "natural"?}\label{zoza}
According to \eref{yop2}, we have 
\begin{equation}\label{pezo}
`E_p(X_x)=`E_p\l(B^{x}(p)\prod_{c : \textrm{ children of }x} (1-X_{c})\r)
\end{equation}
where $B^x(p)$ is the Bernoulli$(p)$ random variable $\1_{C_x=a}$ (recall that all the random variables $B^x(p)$ are independent). 
We now expand the right hand side of \eref{pezo}, conditioning by ${\A}^{\{x\}}$~:
\[`E_p(X_x)=`E_p\l(`E_p\left(B^{x}(p)\prod_{c : \textrm{ children of }x} (1-X_{c}) \Big| {\A}^{\{x\}}\right)\r).\]
Assume now that ${\A}^{\{x\}}$ is known (and finite).
By construction, $B^u(p)=1$ for any $u$ in ${\A}^{\{x\}}$ and $B^u(p)=0$ for any $u \in {\cal P}({\A}^{\{x\}})$· 
We now replace each of the $X_c$ by a product involving its children $B^{c}(p)\,\prod_{c' : \textrm{ children of }c} (1-X_{c'})$, recursively; we stop this expansion when the children of a node are all perimeter sites, since in this case all the $X_{c'}$ equal 0. \par
When one expands this formula (there are no more variables $X_c$ since they have been replaced by those concerning their children, until the perimeter of ${\A}^{\{x\}}$ has been reached), one builds a tree~:  we have $\prod_{c': c'\textrm{ children of }c}(1-X_{c'})=\sum_{C \subset \{\textrm{ children of }c\}} \prod_{c'\in C}-X_{c'}$ corresponding in the tree like decomposition to the choice of some 1 and some $-X_{c'}$. The ``1'' in the decomposition are the leaves of this tree, and when $C$ is chosen, $|C|$ is the degree of an internal node in this tree. 
One then sees that the final value of $X_x$ knowing ${\A}^{\{x\}}$, is
$\sum_{\bt :\pi(\bt) \prec A^{\{x\}}} (-1)^{|\bt| +1}$ in other words $\Delta_x({\A}^{\{x\}})$. 
Indeed, the terms $(-1)^{|\bt| +1}$ counts the number of times $-X$ has been chosen in the expansion~: this number equals the total number of children in the tree, $|\bt|-1$. 
Hence, for $p<p_{crit}^{(x)}$,
\[`E_p(X_x)=\sum_{A}p^{A}(1-p)^{|{\cal P}(A)|}\sum_{\bt :\pi(\bt) \prec A} (-1)^{|\bt| +1}=`E_p(\Delta_x({\A}^{\{x\}}))\]
On the other hand for $p<R_x\leq p_{crit}^{(x)}$,
\[-\bG_x(-p)=\sum_{A}p^{A}(1-p)^{|{\cal P}(A)|}\sum_{A' \prec A} (-1)^{|A|+1}=`E_p(D_x({\A}^{\{x\}})).\]
Then, $D(A)=\Delta(A)=\chi(A)$ is indeed an explication at the level of objects of $-\bG_x(-p)=`E_p(X_x)$.

\subsection{Directed animal with compact sources}

The previous subsections deal almost only with DA with sources of cardinality 1. But most of the results stated there admit some generalizations to DA with larger sources.
The main results that we want to state  is the following generalization of Theorem \ref{fond1}.
\begin{theo}\label{fond4}
Let $G=(V,E)$ be an agreeable graph, $k$ be a positive integer, $S=(s_1,\dots,s_k)$ be a free finite subset of $V$ and $p$ smaller than the radius of convergence of ${\bf G}_S^G$. Under $`P_p$, the random variables $(X_x)_{x\in S}$ are a.s. well defined by \eref{yop2}, and we have
\begin{equation}\label{fonde}
`E_p\l(\prod_{x\in S} X_x\r)=`P_p(X_x=1,x\in S)= (-1)^{|S|}{\bf G}_S^G(-p).
\end{equation}
\end{theo}
Observe that 
\begin{equation}\label{niopse}
`P_p(X_x=1,x\in S)=\sum_{A, A\in {\cal A}(S,G), |A|<+\infty,X_x(A)=1, x\in S}p^{|A|}(1-p)^{|{\cal P}(A)|}
\end{equation}
for any $p\in[0,p_{crit}^S)$ is the sum on  DA occupied on their source $S$.\par
{\bf Proof .}
The first equality is clear, and to prove the second one a generalization of the proof of Theorem \ref{zaea} is needed; we skip the details and give only the main lines of this generalization.\par

First, following the lines of Proposition \ref{etap}, we have
\ben\label{rzz}
-\G^G_{S}(-p)&=&
`E_p\l(D_S({\cal A}^{S})\r)
\een
where for any finite DA $A$,
\begin{equation} \label{diff2}
D_S(A)= \1_{{\cal S}(A)=S}\sum_{A'~:~ A' \prec A, {\cal S}(A')=S} (-1)^{|A'|+1}.
\end{equation}
Let ${\bf f}=(\bt_1,\dots,\bt_{k})\in{\cal T}^k$ be a forest with $k$ trees.  We say that ${\bf f}$ is embeddable in $G$ at $S$ if for any $i\in\{1,\dots,k\}$, $\bt_i$ is embeddable in $G$ at $s_i$. We set
\[\pi_S({\bf f})=\bigcup_{i\in\{1,k\}}\pi_{s_i}(\bt_i)\]
the union of the embeddings of the $ \bt_i$'s on the $s_i's$.
For any forest ${\bf f}$ embeddable in $G$ at $S$, $\pi_S({\bf f})$ is a DA on $G$ with source $S$.
For any DA $A$ with source $S$ in $G$, we set 
\[Y_S(A):=\sum_{{\bf f}=(\bt_1,\dots,\bt_k)~:~{\bf f}\in{\cal T}^k, \pi_S({\bf f})\prec A} (-1)^{|\bt_1|+\dots+|\bt_k|+k},\]
where the sum is taken on the set of forests embeddable in $G$, and
\begin{equation}\label{Delta2}
\Delta_S(A):=\1_{{\cal S}(A)=S}\sum_{A'~:~A'\prec A,{\cal S}(A')=S} Y_S(A').
\end{equation}
A simple analysis yields $Y_S(A):=\prod_{i=1}^k \Delta_{s_i}(A_i)$ 
where $A_i$ is the maximal sub-DA of $A$ with source $s_i$. 
For $A=S$, we have $Y_S(A)=1$. Adapting the proof of Lemma \ref{a} to forests, we get 
\[Y_S(A)=(-1)^{|A|-|S|}.\]
By \eref{diff2} for any $A$ with source $S$, 
\begin{equation}\label{relaze}
D_S(A)=(-1)^{1+|S|}\Delta_S(A).
\end{equation}
It remains to relate $\Delta_S(A)$ with $\prod_{s\in S}X_s$.
By definition of the gas model of type 1, we have
\be
\{`o, X_s(`o)=1, s\in S\} &=& \l\{`o, C_s(`o)=a, s\in S, X_y(`o)=0, y\in S_1\r\}\\
 &=&\l\{`o,{\cal S}(A^{S}(`o))=S, \prod_{y\in S_1}\l(1-\chi_y(A^{(y)}(`o))\r)=1\r\}
\ee
where $S_1$ is those vertices in $G$, children of elements of $S$. 
Write
\be
\{`o,\Delta_S(A^{S}(`o))=1\}&=&\{`o, C_s(`o)=a, s\in S,\prod_{s_i\in S}\Delta_{s_i}(A_i^{(s_i)}(`o))=1\}
\ee
Now, consider $S_1(i)$ the children of $s_i$ in $G$.  By the proof of Theorem \ref{zaea},
\begin{equation}\label{pr}
\Delta_{s_i}(A_i^{(s_i)}(`o))=\1_{\{`o, C_{s_i}(`o)=a\}} \prod_{y\in S_1(i)} (1-\Delta_{y}(A^{y}(`o)))
\end{equation}
and then 
\[\prod_{s_i\in S}\Delta_{s_i}(A_i^{(s_i)}(`o))=\prod_{s_i\in S}\l(\1_{\{`o, C_{s_i}(`o)=a\}} \prod_{y\in S_1(i)} (1-\Delta_{y}(A^{y}(`o)))\r).\]
Now, notice that a node $y$ in $S_1$ may be in several $S_1(i)$'s, but since $1-\Delta_{y}(A^{y}(`o))\in\{0,1\}$ the repetition of this factor in the product \eref{pr} may be simplified; we then obtain
\[\prod_{s_i\in S}\Delta_{s_i}(A_i^{(s_i)}(`o))=\l(\1_{\{`o, C_{s}(`o)=a,\forall s\in S\}} \prod_{y\in S_1} (1-\Delta_{y}(A^{y}(`o)))\r).\]
Using \eref{dix}, this is equal to $\1_{\{`o, C_{s}(`o)=a,\forall s\in S\}} \prod_{y\in S_1}(1-\chi_y(A^{(y)}(`o))$. Finally we get
\[\{`o, X_s(`o)=1, s\in S\}=\{`o,\Delta_S(A^{S}(`o))=1\}\]
and this, together with \eref{relaze} and \eref{rzz} allow to conclude.
$~\Box$

\section{Gas model of type 1 and DA on $\mathbb{Z}^2$}
\label{Z2}
\subsection{New derivation of the GF of DA counted according to the area} 
In this part we work on the square lattice $\Sq=\mathbb{Z}^2$ where as said in the Introduction,  $\Sq$ is viewed as a directed graph where each vertices $(x,y)$ has children $(x,y+1)$ and $(x+1,y+1)$. We represent this lattice as on Fig. \ref{feffe}.
The real number $p_{crit}^{\Sq}:= p_{crit}^{\{x\}}$ does not depend on $x$ by symmetry, and for any finite source $S$, $p_{crit}^{S}=p_{crit}^{\Sq}$. According to Proposition \ref{yopop}, $p_{crit}^{\Sq}\geq 1/2$ and by Proposition \ref{tru}, for $p<p_{crit}^{\Sq}$ the gas occupation is $`P_p$ a.s. defined everywhere. The radius $R_x$ can also be shown easily to be positive using a simple surjection from the set of DA with $n$ cells in $\Sq$ into the set of trees having only internal nodes with total degree 2 and 3 and $n$ nodes.\par

The set $L_i:=\l\{(x,i-x),x\in\mathbb{Z}\r\}$ is called the $i$th horizontal line and we denote by  $Z_i$ its gas occupation~: 
\[Z_i(x)= X_{(x,i-x)},~~~\textrm{ for any }x\in\mathbb{Z}.\]
By construction $Z_i$ is a process indexed by $\mathbb{Z}$ taking its values in $\{0,1\}$. It is invariant in the strong sense~: for any $k\in \mathbb{Z}$ the process $(Z_i(x))_{x\in \mathbb{Z}}$ and  $(Z_i(x+k))_{x\in \mathbb{Z}}$ have the same distribution. By construction also the processes $Z_i$ and $Z_{i+1}$ have the same distribution. Since $Z_{i+1}$ is built from $Z_i$ via the gas evolution \eref{yop2}, for any $i\in\mathbb{Z}$, $Z_i$ and $Z_{i+1}$ are related by~:
\begin{equation}\label{rezr}
Z_i(k)=B_k^i(p) \, (1-Z_{i+1}(k))\,(1-Z_{i+1}(k+1))
\end{equation}
where the $B_k^i(p)$ are independent and identically distributed (i.i.d.) Bernoulli$(p)$ random variables, independent of $Z_{i+1}$.

The existence (and the construction) of a solution is guaranteed when $p\in[0,p_{crit}^{\Sq})$ by Proposition \ref{tru}, and a solution, say $Z^\star$ is characterized by the equation \eref{niopse}~: let $S=\{(x,-x), x\in I\}$ be included in $L_0$, with $I$ finite. We have
\[`P(X_s=1, s\in S)=`P(Z^\star(l)=1,l\in I)=\sum_{A,A\in {\cal A}(S),\chi_s(A)=1,s\in S}p^{|A|}(1-p)^{|{\cal P}(A)|}.\]
By inclusion-exclusion, this determines the finite dimensional distribution of $Z^{\star}$. 
\begin{pro}\label{uni}
 For all $p\in[0,p_{crit}^{\Sq})$, there exists a unique law $\mu$ of processes taking their values in $\{0,1\}$ solution of \eref{rezr} (in other words, there exists a unique law $\mu$ such that if the process $Z_{i+1}$ is $\mu$ distributed then so do $Z_i$).
\end{pro}
\proof The uniqueness of the solution will be proved to be a consequence of Theorem \ref{fond4}~: assume that $Z'_0$ has a certain law $\mu'$ which is solution of \eref{rezr}. One then builds an infinite sequence $(Z'_i)_{i\leq -1}$ such that for any $Z'_i$ is obtained from $Z'_{i+1}$ by
\begin{equation}\label{rezrr}
Z'_i(k):=B_k^{(i)}(p) \, (1-Z'_{i+1}(k))\,(1-Z'_{i+1}(k+1)) \textrm{ for }i\leq -1
\end{equation}
where $B_k^{(i)}(p), i\leq -1, k\in\mathbb{Z}$ is an array of i.d.d. Bernoulli$(p)$ random variables. \par
\noindent Consider $I$ a finite subset of $\mathbb{Z}$. We will examine the probability $`P(Z'_i(l)=1, l\in I)$ and show that it converges when $i\to+\infty$ to the distribution of  $`P(Z^\star_i(l)=1, l\in I)$. Since $`P(Z'_i(l)=1, l\in I)$ does not depend on $i$ this entails that the finite distribution of $Z'$ and those of $Z^\star$ are equals. This result being valid for any finite set $I$, this implies that $Z'$ has the same law as $Z^\star.$ \par

Let $I$ be fixed, and consider $S(i)=\{(x,i-x),x\in I\}$ the subset of the $i$th line $L_i$, which is simply obtained by the translation of $S:=S(0)$. Under $`P_p^\Sq$, $\A^S$ is a.s. finite; this implies that for any $`e>0$ there exists $n_{`e}$ such that $`P_p(|\A^S|\geq n_{`e})<`e$.\par
Consider now $S(-n_{`e})$, the subset of $L_{-n_{`e}}$. We have  $`P_p(|\A^{S(-n_{`e})|}\geq n_{`e})<`e$, by invariance by translation. When $|\A^{S(-n_{`e})}|< n_{`e}$, $\left(\A^{S(-n_{`e})}\cup {\cal P}(\A^{S(-n_{`e})})\right)\cap L_0=\emptyset$. In this case, the values of $(X_s)_{s\in S(-n_{`e})}$ does not depends on the value of $X$ on $L_0$ (and then neither on its distribution on this line), since as said in Section \ref{cgm}, the gas occupation on a subset $S$ is a deterministic function of $\A^S$ and its perimeter sites (and then does not depend on the other sites). \par
The formula \eref{rezrr} defining $Z'_{-n}$ is the same as that defining the gas model of type 1, and then, if $\A^{S(-n_{`e})}$ 
satisfies $|\A^{S(-n_{`e})}|< n_{`e}$, then $(Z'_{-n_{`e}}(l),l\in I)$ does not depend on $Z'_0$. 
Therefore, for any $n\geq n_{`e}$, we have
\[\left|`P_p(Z'_{-n}(l)=1, l\in I)-`P_p(Z^\star(l)=1, l\in I)\right|<`e \] 
and then $`P_p(Z'_{-n}(l)=1, l\in I)\to`P_p(Z^\star(l)=1, l\in I).~\Box$
\begin{figure}[htbp]
\psfrag{z0}{$Z'_0 \sim \mu'$  }
\psfrag{z1}{$Z'_{-1}$}
\psfrag{z2}{$Z'_{-2}$}
\psfrag{z3}{$Z'_{-3}$}
\psfrag{z4}{$Z'_{-4}$
}\psfrag{sources1}{Observed sites does not depend on $Z'_0$}
\centerline{\includegraphics[height=4cm]{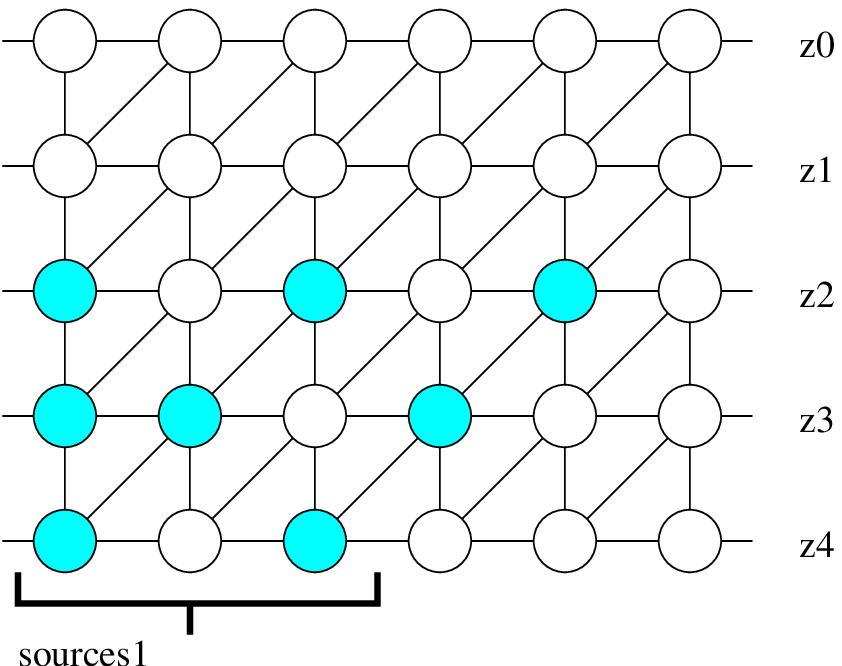}}
\captionn{\label{diffs} The colored sites are the sites where the Bernoulli$(p)$ coloring $=1$. The gas occupation on the three first cells on the last row depends only on the DA  (with these three cells as over-sources). The height of this DA is so small that the gas occupation of these cells does no depend on $Z'_0$. }
\end{figure}

\begin{com}The sequence $(Z_i)$ forms  (``a vertical'') Markov chain with infinite number of states. 
If one wants to determine the value of $Z_i(x)$ using the Bernoulli tossing allowing to realize this Markov chain, only the DA of calculus ${\A}^{(x,i-x)}$ is needed. This is reminiscent of the Propp \& Wilson coupling from the past technique to simulate a Markov chain under the stationary distribution~: only the last transitions of the Markov chain sufficient to determine the current state are needed. 
\end{com}

We now determine the distribution $\mu$  of the processes $Z_i$  and the density of the gas model of type 1.

\begin{figure}[htbp]
\psfrag{b1}{$B^\circ_1$}\psfrag{b3}{$B^\circ_2$}\psfrag{b5}{$B^\circ_3$}\psfrag{b2}{$B^\bullet_1$}\psfrag{b4}{$B^\bullet_2$}\psfrag{t1}{$\tau^{\bullet}_1$}\psfrag{0}{0}\psfrag{1}{1}
\centerline{\includegraphics[height=1.5cm]{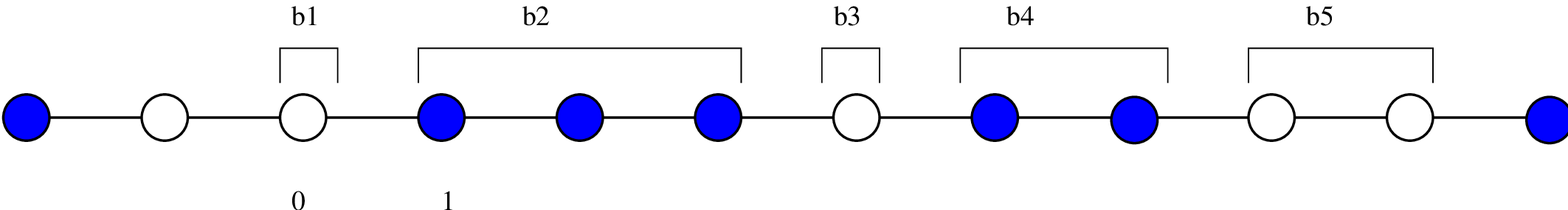}}
\captionn{\label{tb} On this example, the black cells are occupied and corresponds to the places where $Z=1$, and the whites one are empty.}
\end{figure}

Let $Z:=Z_0=(Z_0(x))_{x\in \mathbb{Z}}$ be the gas process on the line. 
We denote by $(B_i^\bullet)_{i\geq 1}$ (resp. $(B_i^\circ)_{i\geq 1}$) the successive sizes of the blocks of consecutive occupied (empty) positions at the right of zero. If a block contains some negative positions, only the part at the right of zero is counted as  on Fig. \ref{tb}.

 A random variable $G$ is said to follow the geometric distribution with parameter $\alpha$, we denote $ G \sim{\cal G} (\alpha)$, when $`P(G=k)=\alpha(1-\alpha)^{k-1} \textrm{ for any } k\geq 1.$
\begin{theo}\label{zoa}
 Let $p\in(0,p_{crit})$. The density of the gas model is given by
\[`P_p(Z(0)=1)=\frac{1/\alpha_{\bullet}}{1/\alpha_\circ+1/\alpha_\bullet},\]
where\[\alpha_\bullet=\frac{-1+p+\sqrt{1+2p-3p^2}}{2p}~~~~\textrm{ and }~~~~\alpha_\circ=\frac{-1+p+\sqrt{1+2p-3p^2}}{1+p+\sqrt{1+2p-3p^2}}.\] 
Under $`P_p$ the random variables $B_{i}^{\bullet}$ and $B_i^{\circ}$ for $i\geq 1$ are independent and independent of $Z(0)$, and $B_i^{\bullet}\sim {\cal G}(\alpha_\bullet)$ and  $B_i^{\circ}\sim {\cal G}(\alpha_\circ)$.
\end{theo}
Hence, for any $b^\circ_i,b^\bullet_i\geq 1$
\[`P(Z_0(0)=x,B^\circ_i=b^\circ_i,B^\bullet_i=b^\bullet_i,i\in\{1,\dots,k\})=\frac{1/\alpha(x)}{1/\alpha_\bullet+1/\alpha_\circ}\prod_{i=1}^k `P(G^\circ=b^\circ_i)`P(G^\bullet=b^\bullet_i)\]
where $G^\bullet\sim {\cal G}(\alpha_\bullet)$ and $G^\circ\sim {\cal G}(\alpha_\circ)$, and $\alpha(1)=\alpha_\bullet$ and $\alpha(0)=\alpha_\circ$.\par
 The law $\mu$ of $Z$ is characterized by the properties given in this Theorem and the invariance by translation, since this characterizes the finite dimensional distributions of $\mu$.  
We have no combinatorial explanation to the properties of the blocks, but these distributions may be guessed from the result of Bousquet-M\'elou cited at the beginning of Section \ref{eraz}.
\begin{com}\label{ezsd}
$\bullet$ The proof that we provide below shows that the distribution described in Theorem \ref{zoa} is a solution of equation \eref{rezr} not only for $p\in(0,p_{crit})$ but also for any $p\in(0,1)$. For $p\in [p_{crit},1)$ we didn't find an argument to show that there is a unique solution to \eref{rezr}. In the case $p=1$, the solution is not unique since an alternating sequence of blocks, the occupied block having size $k$ and the empty one, $k+1$ is conserved by the gas evolution for any $k\geq 1$. For $p=0$, the process $Z$ null a.s. is the unique solution.\\
$\bullet$ One may also view $Z$ as built as  juxtaposing alternatively occupied and empty blocks with the prescribed distributions, starting from $-\infty$. Viewed like this, the density can be computed using a simple renewal argument or by the law of large number.\\
$\bullet$ Assume that $Z_1=(Z_1(x))_{x\in \mathbb{Z}}$ and $Z_0=(Z_0(x))_{x\in \mathbb{Z}}$ have the property announced in the Theorem, but 
that $\alpha_\bullet$ and $\alpha_\circ$ are unknown.
To compute  $\alpha_\bullet$ and $\alpha_\circ$ we derive and solve the two following equations  \begin{equation}\label{yoop}(1-\alpha_{\circ})p=1-\alpha_{\bullet} \textrm{ and }\alpha_{\circ}=\alpha_{\bullet}(1-\alpha_{\circ})p.
\end{equation} 
Indeed, a block of $1$ in $Z_0$  is placed ``under'' a block of $0$ in $Z_1$. Once a block of 1 is begun in $Z_0$, there is an additional 1 at the right of this block with probability $(1-\alpha_\circ)p$ (an additional 0 is needed in $Z_1$ and the right Bernoulli tossing is needed). This gives the first equation, since in $Z_1$ an additional 1 occurs with probability $1-\alpha_\bullet$. On the other hand,
since $\{x, Z_0(x)=1\}= \b\{x, C_x:=a, Z_1(x)=0, Z_1(x+1)=0\b\}$, the probability $\frac{1/\alpha_\bullet}{1/\alpha_\bullet+1/\alpha_\circ}=`P_p(Z_0(x)=1)=p`P(Z_1(x)=0, Z_1(x+1)=0)$ equals to $p`P(Z_1(x+1)=0 |Z_1(x)=0 )`P(Z_1(x)=0)=p \frac{1/\alpha_\circ}{1/\alpha_\bullet+1/\alpha_\circ}(1-\alpha_0)$.  
This give the second equation (that can also be obtained by a renewal type argument without using the density).
\end{com}

\subsection{Consequences of Theorem \ref{zoa}}
It is worth mentioning that Theorem \ref{zoa} is in fact equivalent to the following representation, allowing also to compute the gas density. 
\begin{pro}\label{tyr}
Under $`P_p$, the process $(Z(x))_{x\in\mathbb{Z}}$ is a Markov chain with transition matrix 
\begin{equation}\label{M}
M=\left(\begin{array}{cc}
`P(Z(1)=1 | Z(0)=1)& `P(Z(1)=0 | Z(0)=1)\\ `P(Z(1)=1 | Z(0)=0)& `P(Z(1)=0 | Z(0)=0)
\end{array}\right)=\left(
\begin{array}{cc} 
1-\alpha_{\bullet} & \alpha_\bullet\\
\alpha_\circ & 1-\alpha_\circ
\end{array}\right)
\end{equation}
under the stationary distribution.
\end{pro}
One may compute from this Proposition and Theorem \ref{fond4} some results about the GF of DA with general sources on the square lattice.
\begin{pro}\label{geos} Let $S=\{s_1,\dots,s_k\}$ where $s_i=(x_i,-x_i)$ be some points on the principal line $L_0$, such that $d_i:=x_{i+1}-x_i$ for $i\in\{2,\dots,k\}$ are positive integers. The GF of DA on the square lattice with source $S$ is given by
\[{\bf G}_{S}^\Sq(-p)=(-1)^{|S|}\frac{1/\alpha_\bullet}{1/\alpha_\bullet+1/\alpha_\circ}\prod_{i=1}^{k-1} \frac{\alpha_\bullet(1-\alpha_\bullet-\alpha_\circ)^{d_i}+\alpha_\circ}{\alpha_\bullet+\alpha_\circ}\]
\end{pro}
\proof This is a simple consequence of the Markovianity of $Z$ and of Theorem \ref{fond4}. Using a diagonalization of $M$ (given in \eref{M}), one gets
$`P(Z_k=1 | Z_0=1)= \l(M^k\r)_{1,1}=\frac{\alpha_\bullet(1-\alpha_\bullet-\alpha_\circ)^{k}+\alpha_\circ}{\alpha_\bullet+\alpha_\circ}$.~$\Box$\medskip

\begin{figure}[htbp]
\psfrag{x1}{$d_1$}
\psfrag{x5}{$d_5$}
\psfrag{s1}{$s_1$}
\psfrag{s5}{$s_5$}
\psfrag{x2}{$d_2$}
\psfrag{x3}{$d_3$}\psfrag{s2}{$s_2$}\psfrag{s3}{$s_3$}
\psfrag{d3}{$d_3$}
\psfrag{x4}{$d_4$}
\centerline{\includegraphics[height=2.5cm]{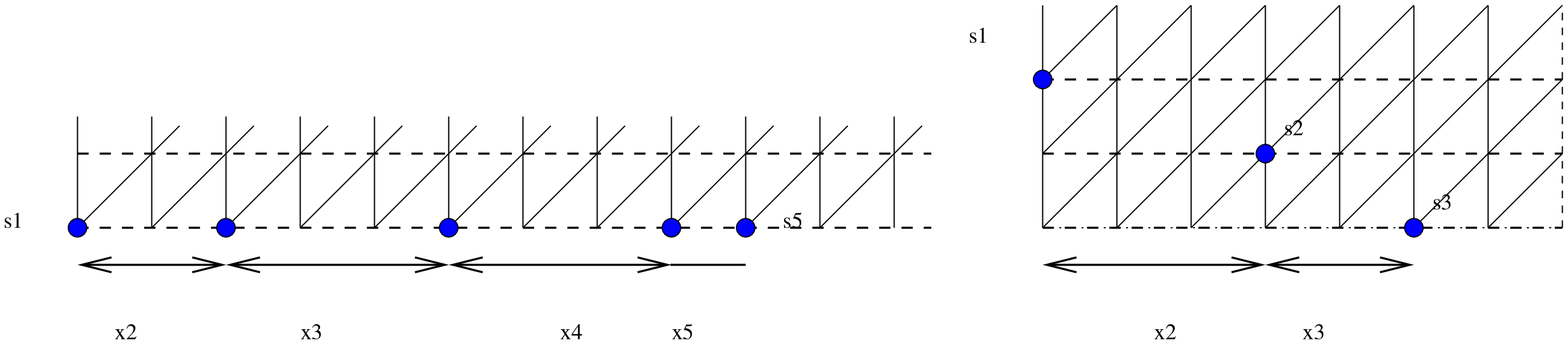}}
\captionn{\label{source}On the first picture, the sources considered in Proposition \ref{geos}. The second picture illustrates the sources considered in Proposition \ref{geos2}.}
\end{figure}

For example, if $S_n:=\{(i,-i), i=1,\dots,n\}$ then 
${\bf G}_{S_n}^\Sq(-p)=\frac{1/\alpha_{\bullet}}{1/\alpha_{\bullet}+1/\alpha_{\circ}}(1-\alpha^{\bullet})^{n-1}(-1)^{n}.$
Then, the GF of DA on the square lattice with compact sources satisfies
\[\sum_{n\geq 1} {\bf G}_{S_n}^\Sq(-p)=\frac{1/\alpha_{\bullet}}{1/\alpha_{\bullet}+1/\alpha_{\circ}}\sum_{n\geq 1} (1-\alpha^{\bullet})^{n-1}(-1)^{n}=\frac{-p}{1+3p}. 
\]
A combinatorial explanation of this formula is given in Gouyou-Beauchamps \& Viennot \cite{GBV}. \medskip
One may also compute the GF of DA having their sources on different lines. For example~:
\begin{pro}\label{geos2} 
Let $S=\{s_1,\dots,s_k\}$ where $s_i=(x_i,-y_i)$ be some points of $\Sq$ such that $y_i=x_i+i$ (see Figure \ref{source} $(ii)$) and $d_i:=x_{i+1}-x_i$ for $i\in\{2,\dots,k\}$ are positive integers. The GF of DA on the square lattice with source $S$ is given by
\[{\bf G}_{S}^\Sq(-p)=(-1)^{|S|}\frac{1/\alpha_\bullet}{1/\alpha_\bullet+\alpha_\circ}\prod_{i=1}^{k-1}{\alpha_\circ} \frac{(1-\alpha_\bullet-\alpha_\circ)^{d_i}-1}{\alpha_\bullet+\alpha_\circ}\]
\end{pro}
See Figure \ref{source} for an illustration of the considered sources.\par
\proof We use again the Markovian properties of $Z$ but this time on several lines in the same time. We want to compute $`P(X_s=1,s\in S)$. Notice that the cells of $S$ are not ancestors from each others. For any $s\in S$, denote by $c_1(s)$ and $c_2(s)$ the two children of $s$ in $\Sq$ ($c_1$ is at the left of $c_2$); we have 
\[`P(X_s=1,s\in S)=`P(X_s=1,X_{c_1(s)}=X_{c_2(s)}=0,s\in S)\] 
since if $X_s=1$ then  a.s. $X_{c_1(s)}=X_{c_2(s)}=0$. We condition on the gas occupation on $s_1$, $c_1(s_2)$ and $c_2(s_2)$~:
\be
`P(X_s=1,s\in S)&=&`P(X_s=1, s\in S\setminus\{s_1\} | X_{s_1}=1,X_{c_1(s_2)}=0,X_{c_2(s_2)}=0)\\
&&~\times`P( X_{s_1}=1,X_{c_1(s_2)}=0,X_{c_2(s_2)}=0).\ee
By Markovianity we have $`P( X_{s_1}=1,X_{c_1(s_2)}=0,X_{c_2(s_2)}=0)=\frac{1/\alpha_\bullet}{1/\alpha_\bullet+1/\alpha_\circ} (M^{d_1})_{1,2}M_{2,2}$
and 
\be
`P(X_s=1, s\in S\setminus\{s_1\} | X_{s_1}=1,X_{c_1(s_2)}=X_{c_2(s_2)}=0)&=&`P(X_s=1, s\in S\setminus\{s_1\}|X_{c_1(s_2)}=X_{c_2(s_2)}=0 )\\
&=&\frac{`P(X_s=1, s\in S\setminus\{s_1\})}{`P(X_{c_1(s_2)}=X_{c_2(s_2)}=0)}
\ee
Thus, $`P(X_s=1,s\in S)=`P(X_s=1, s\in S\setminus\{s_1\}) \dis\frac{\alpha_\circ}{\alpha_\bullet}(M^{d_1})_{1,2}$. Using that $(M^k)_{1,2}=\frac{\alpha_\bullet(1-\alpha_\bullet-\alpha_\circ)^{k}-\alpha_\bullet}{\alpha_\bullet+\alpha_\circ}$ we get the result.
~$\Box$

\begin{remi}\label{reme} Let $S=\{s_1,\dots,s_k\}$ be a free subset of $\Sq$. The serie $\bG_S$ may be computed (theoretically) using the properties of the gas model, since it suffices to compute $`P(X_s=1,s\in S)$. This can be done by writing a finite sum involving the occupations of the cells that are above the $s_i$'s and under the line $L_j$ with the largest index such that $L_j\cap S$ is not empty (inside the region surrounded by the dotted lines in Fig. \ref{soce}), and using the Markovianity of the gas occupation on $L_j$.
\begin{figure}[htbp]

\centerline{\includegraphics[height=2cm]{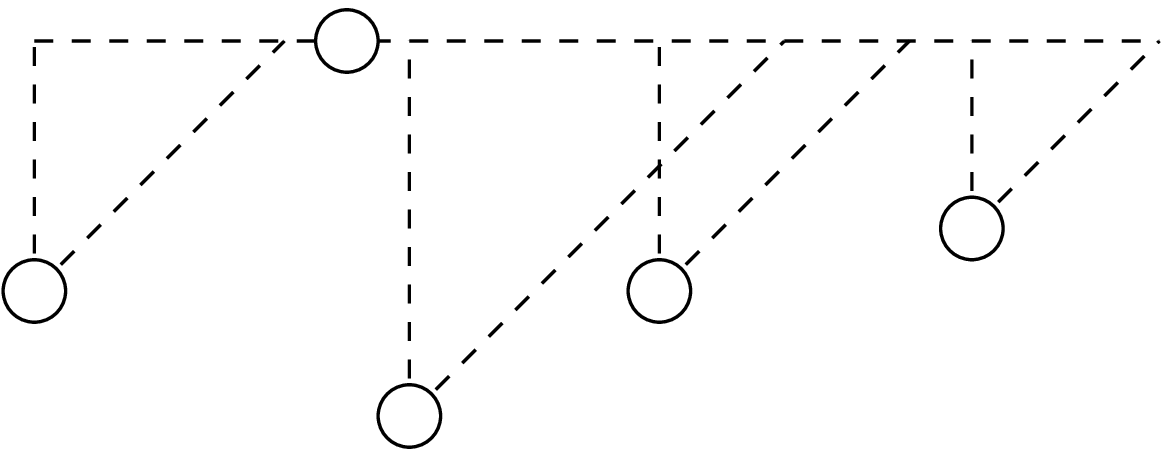}}
\captionn{\label{soce}The points stands for the points of $S$; the dotted lines surround the (finite) set of cells where the summation has to be done.}
\end{figure}

\end{remi}

The Markovianity of $Z$ implies some (weighted) properties of DA taking into account the occupation or not of their sources; here, by weighted we mean that here we deal with probabilities, and then each DA $A$ has weight $p^{|A|}(1-p)^{|{\cal P}(A)|}$ which depends on its area and perimeter. Let $z:=(z_0,\dots,z_k)$ be fixed in $\{0,1\}$.
The probability of the event $(Z(0),\dots,Z(k)) = z$ depends only on the four quantities $n^{\bullet\bullet}_k(z) = \#\{i, i\in\{0\ldots k-1\} , (z_i,z_{i+1}) = (1,1)\}$, $n^{\bullet\circ}_k(z)$, $n^{\circ\bullet}_k(z)$ and $n^{\circ\circ}_k(z)$, defined accordingly. 
Hence $`P((Z(0),\dots,Z(5))=z)=`P((Z(0),\dots,Z(5))=z')$ for $z=(0,0,1,1,0,0)$ and $z'=(0,1,1,0,0,0)$. In terms of DA with horizontal sources $S=\{(i,-i)\}_{i=0\ldots 5}= \{v_i\}_{i=0\ldots 5}$ it leads to 
\[\sum_{A, \chi_{v_i}(A) = z_i, i\in\{0,\ldots,5\}} p^{|A|}(1-p)^{|\mathcal{P}(A)|} = \sum_{A, \chi_{v_i}(A) = z'_i,i\in\{0,\ldots,5\}} p^{|A|}(1-p)^{|\mathcal{P}(A)|}  \]
A bijection between DA with the considered occupation of their sources and preserving moreover the area and the perimeter would be a direct explanation of this identity, but such a map does not exist since the enumeration of these two sets of DA according to the area and the perimeter are distinct for DA of area $5$ :
\[ \sum_{A, |A| \leq 5,\chi_{v_i}(A) = z_i, i\in\{0,\ldots,5\}}  t^{|A|}u^{|\mathcal{P}(A)|} = u^7t^2+(2u^8+10u^9)t^4+(2u^8+7u^9)t^5 \]   
and
\[  \sum_{A, |A| \leq 5,\chi_{v_i}(A) = z'_i, i\in\{0,\ldots,5\}}   t^{|A|}u^{|\mathcal{P}(A)|} =  u^7t^2+(2u^8+10u^9)t^4+(1u^8+8u^9)t^5.\]

\subsection{Proof of Theorem \ref{zoa}}
\label{eraz}

 \rm One may propose several proofs of Theorem \ref{zoa} ~: following the proof of Bousquet-Mélou \cite{MI1}, one may first work on the cylinder as done in Section \ref{SQAM}. There the process $Z$ is in some sense horizontally Markovian, and conditioned to come back at its starting point. Theorem \ref{zoa} may be obtained by a passage to the limit on the length of the cylinder (the limit in distribution here is just a convergence of the finite dimensional distribution). \par
More directly, another method consists in proving first Proposition \ref{tyr}~: three steps are needed. At the beginning, suppose that $Z_1$ is a Markov chain with two states $0$ and $1$. Then using that is must be preserved by the gas evolution, there is at most one Markov chain satisfying these constraints (easy). Let $M$ be its transition matrix. It remains to show that this Markov chain is indeed preserved by the gas evolution. This is done via some linear algebra on $M$, showing that if $(Z_1(i),\dots,Z_1({i+k}))$ is a $M$-Markov chain, so do $(Z_0(i),\dots,Z_0({i+k-1}))$ (this is done via some quite laborious computation). Proposition \ref{uni} allows to conclude.
\medskip

The proof of Theorem \ref{zoa} that we propose below may be viewed as a translation of what is said above. But we think that it is in fact different and is more likely to be generalized to other models. Another pleasant reason is the natural appearance of the gas density in the beginning of the considerations (see Comment \ref{ezsd})~: using the other approaches, this is not the case.\medskip

\noindent \bf Proof of Theorem \ref{zoa}. \rm
Let $(b_i^\circ)_{i\in\{1,\dots,k\}}$ and $(b_i^\bullet)_{i\in\{1,\dots,k\}}$ be two fixed vectors in $\{1,2,\dots\}^{k}$.
For $j$ in $\{0,1\}$ we denote the block sizes of  $Z_j$ by $B_{i,j}^\circ$ and $B_{i,j}^\bullet$. 
We assume that $Z_1(0)$ and the block sizes of $Z_1$ have the distribution described in the Theorem, and we prove then that $Z_0(0)$ and the blocks of $Z_0$ have the same distribution. By Proposition \ref{uni} this is sufficient to prove the Theorem.\\

\begin{figure}[htbp]
\psfrag{Z0}{$Z_0$}\psfrag{Z1}{$Z_1$}\psfrag{k}{$k$}\psfrag{k+1}{$k+1$}
\centerline{\includegraphics[width=12 cm]{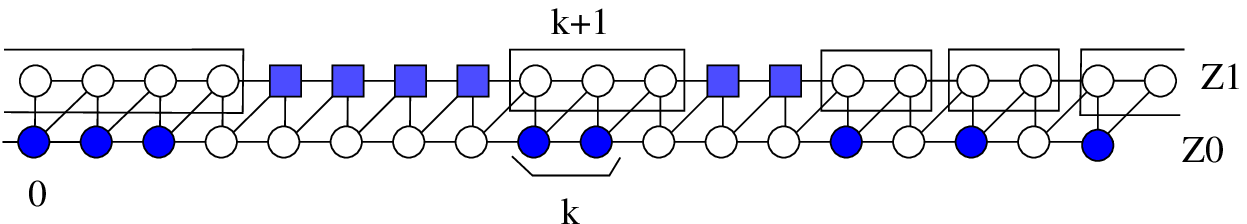}}
\captionn{\label{tbtb}A block of $k$ occupied positions in $Z_0$ determines $k+1$ empty position on $Z_1$ (surrounded by rectangles). The little squares figure out the uncertain cells, that can not be determined just in view of the gas occupation in $Z_0$.}
\end{figure}

To compute 
$`P({Z}_0(0)=1, {B}_{i,0}^\circ=b_i^\circ, {B}_{i,0}^\bullet=b_i^\bullet,i\in\{1,\dots,k\})$, we will sum on all possible block configurations of $Z_1$ on the positions $\{0,\dots,1+\sum_{i=0}^k b_i^\circ+b_i^\bullet\}$. We recall that by the gas model of type 1, the process $Z_0$ is solution of \eref{rezr} (see Fig. \ref{tbtb}). 
Thus, knowing $k$ consecutive occupied cells in $Z_0$ allows to guess $k$ Bernoulli random variables and a block of $k+1$ empty positions in $Z_1$ (that maybe included in a larger empty block). In particular, if $Z_0(0)=1$ then $Z_1(0)=0$.\par

The empty blocks of $Z_0$ are a little bit more tricky to handle~: for each empty block in $Z_0$ with size $b_i^\circ$, there are $b_i^{\circ}-1$ cells in $Z_1$ that can be occupied or empty (see Fig. \ref{tbtb}).  We call these cells \it uncertain cells\rm. \par
Notice also that $\{{Z}_0(0)=1, {B}_{i,0}^\circ=b_i^\circ, {B}_{i,0}^\bullet=b_i^\bullet,i\in\{1,\dots,k\}\}=\{{Z}_0(0)=1, {B}_{i,0}^\circ=b_i^\circ, {B}_{i,0}^\bullet=b_i^\bullet,i\in\{1,\dots,k\},{Z}_0(\sum_{i=1}^kb_i^\bullet+b_i^\circ)=1\}$ since knowing the size of the $2k$ first blocks and the values of $Z_0(0)$ implies the occupation of the cell that follows.
Hence, we have
\[`P(Z_0(0)=1, B_{i,0}^\circ=b_i^\circ,B_{i,0}^\bullet=b_i^\bullet, i\in\{1,\dots,k\})=\frac{1/\alpha_\circ}{1/\alpha_\circ+1/\alpha_\bullet}[(1-\alpha_{\circ})p]^{\sum_{i=1}^k b_i^\bullet}\l(\prod_{i=1}^k S_{b_i^\circ}\r)(1-\alpha_0)p.\]
The first factor $\frac{1/\alpha_\circ}{1/\alpha_\circ+1/\alpha_\bullet}$ on the right hand side comes from $Z_1(0)=0$. The other term are computed successively, as transition probabilities~: \\
-- the term  $(1-\alpha_{\circ})p$ corresponds to the creation of an occupied cell $x$ in $Z_0$ knowing that $x$ is empty in $Z_1$~: 
\be
`P(Z_0(x)=1 | Z_1(x)=0)&=&\frac{`P(Z_0(x)=1,Z_1(x)=0)}{`P(Z_1(x)=0)}=\frac{`P(Z_0(x)=1,Z_1(x)=0,Z_1(x+1)=0)}{`P(Z_1(x)=0))}\\
&=&`P(Z_0(x)=1|Z_1(x)=0,Z_1(x+1)=0)`P(Z_1(x+1)=0|Z_1(x+1)=0)\\
&=&p(1-\alpha_{\circ}).
\ee 
In other words knowing that $x$ in empty in $Z_1$ it is occupied in $Z_0$ if there is a transition empty-empty in $Z_1$ between $x$ and $x+1$ and a favorable coloring $C_x=a$ of $x$ (proba. $p$). There are ${\sum_{i=1}^k b_i^\bullet}$ such terms.\\
-- The last term $(1-\alpha_0)p$ comes from the creation of an occupied cell after these blocks in $Z_0$.\\
-- The factor $S_{b_i^\circ}$ corresponds to the contributions of the uncertain cells in $Z_1$ above the $b_i^\circ$ empty cells of $i$th empty considered block of $Z_0$~: more precisely 
\ben\label{lic}
S_{l}=`P(Z_0(1)=0,\dots,Z_0(l)=0,Z_1(l+1)=0 | Z_1(1)=0).\een

Above the empty cells $1,\dots,l$ in $Z_0$ there are $l-1$ incertains cells in $Z_1$, at position $2,\dots,l$. \par
For $l=1$ there is no uncertain cell above and we find $S_1 = (1-\alpha_\circ)(1-p)$. For $k\geq 2$, we write $S_k=S_k^{\bullet}+S_k^\circ$, where $S_k^{\bullet}$ (resp. $S_k^{\circ}$) corresponds to formula \eref{lic}, where in the RHS is added $Z_1(l)=1$ (resp.  $Z_1(l)=0$), in other words where the occupation of the last uncertain cell is specified.
For $k=1$, we take the convention $S_1^{\circ} = S_1$ and $S_1^{\bullet} = 0$. For $k \geq 2$, there is a simple decomposition of $S_k^{\bullet}$ and $S_k^{\circ}$ obtained by removing the last cells of this block:
\[\left\{\begin{array}{ccl}
S_l^\circ&=&S_{l-1}^\circ(1-\alpha_\circ)(1-p)+S^\bullet_{l-1}(1-\alpha_\circ)(1-p)\\
S_l^\bullet&=&\dis S_{l-1}^\circ\frac{\alpha_\circ\alpha_\bullet}{(1-\alpha_\circ)(1-p)}+S_{l-1}^\bullet(1-\alpha_\bullet)
\end{array}\right.\]
Using \eref{yoop}, this finally gives $S_l^\circ=(1-p)(1-\alpha_\circ)S_{l-1}$ and $S_l^{\bullet}=(1-\alpha_\circ)pS_{l-1}$, then $S_l=(1-\alpha_\circ)^{l-1}S_1=(1-\alpha_\circ)^l(1-p)$. Thus  
\[\begin{array}{l}
`P(Z_0(0)=1, B_{i,0}^\circ=b_i^\circ,B_{i,0}^\bullet=b_i^\bullet, i\in\{1,\dots,k\}) \\
\displaystyle \hspace{5cm} = \hspace{1cm} \frac{(1-\alpha_\circ)p/\alpha_\circ}{{1}/{\alpha_\circ}+{1}/{\alpha_\bullet}}[(1-\alpha_{\circ})p]^{\sum_{i=1}^k b_i^\bullet}(1-\alpha_\circ)^{\sum_{i=1}^k {b_i^\circ}}(1-p)^{k}
\end{array}\]
using that $(1-\alpha_\circ)(1-\alpha_\bullet)(1-p)=\alpha_\circ\alpha_\bullet$, and $(1-\alpha_\circ)p/\alpha_\circ=1/\alpha_\bullet$, we get
\be
`P(Z_0(0)=1, B_{i,0}^\circ=b_i^\circ,B_{i,0}^\bullet=b_i^\bullet, i\in\{1,\dots,k\})
&=& \frac{1/\alpha_\bullet}{1/\alpha_\circ+1/\alpha_\bullet}\prod_{i=1}^k`P(G^\bullet=b_i^\bullet)`P(G^\circ=b_i^\circ)
\ee
which is the expected result.~ ~\\
The computation of $`P(Z_0(0)=0, B_{i,0}^\circ=b_i^\circ,B_{i,0}^\bullet=b_i^\bullet, i\in\{1,\dots,k\})$ can be performed along the same lines, except that here the two cases $Z_1(0)=1$ and $Z_1(1)=1$ have to be considered.
$\Box$

\section{Simultaneous construction of DA and gas model of type 2}
\label{gas2}

We construct in this Section a probability space on $G$ in such a way that the density of the gas model equals up to some change of variables the area and perimeter GF for  DA.
This approach is very similar to that of Section \ref{CAG}.
Here a probability measure having two parameters is needed. The idea is again to transfer a quantity  revealing the perimeter size via a DA of calculus. The main difference with Section \ref{CAG} is that here, the revealing quantity is random when it was deterministic in Section \ref{CAG} (it was $\chi(A)$). We endow the perimeter sites with i.i.d Bernoulli$(p)$ random variables, and we will transfer to the source the minimum of those random variables.  Assume that $N$ is random and that you know the law of $M:=\min\{B_1,\dots,B_N\}$ for $B_1,\dots,B_N$ i.i.d. Bernoulli$(p)$ random variables, independent of $N$. It is straightforward that the distribution of $M$ characterizes that of $N$. This is morally what we use here. \par
Consider $G$ an agreeable  graph, and the probability space $\Omega^G_\star=\{a,b,c\}^G$, endowed with the probability measure 
\[`P_{p_a,p_b,p_c}=(p_a\delta_{a}+p_b\delta_b+p_c\delta_c )^{\otimes G},\]
where $p_a,p_b,p_c$ are three non negative parameters summing to 1.
In other words, we have again a random coloring of $G$, where the nodes colors are i.d.d. and follows the following distribution (we write $C^{\star}_x$ the color of $x$). Under $`P_{p_a,p_b,p_c}$, for any $x$,  $C^{\star}_x$ is $a$ with probability $p_a$, $C^{\star}_x$ is $b$ with probability $p_b$ and  $C^{\star}_x$ is $c$ with probability $p_c$.

We define now $\A_\star^{S}(`o)$ in the same way as $\A^{S}(`o)$ is defined in the beginning of Section \ref{COD}~: $\A_\star^{S}(`o)$ has again as set of cells the maximal DA with over-source $S$ whose cells are all $a$-colored; now the perimeter sites of $\A^{S}_\star$ may have the color $b$ or $c$. In order to take into account these colors, we introduce
\[{\cal P}_b(A)=\{x, x\in{\cal P}(A), C^{\star}_x=b\},~~\textrm{ and }~~ {\cal P}_c(A)=\{x, x\in{\cal P}(A), C^{\star}_x=c\},\]
the perimeter sites of $A$ having color $b$ and $c$. The set of DA on $G$ with bi-colored perimeter sites and source $S$ is denoted by ${\cal A}^\star(S,G)$, and we add in the same way a star to the variables already defined in the previous section to make visible that we are working with colored objects.

The following proposition is the analogous, or more precisely a refinement, of Proposition \ref{moch}~:
\begin{pro}Let $G=(V,E)$ be an agreeable graph.\\  $(i)$ Let $B\in {\cal A}^\star(S,G)$ be a finite DA with source (exactly) the free set $S$. We have
\[`P_{p_a,p_b,p_c}({\A}_\star^{S}=B)=p_a^{|B|}p_b^{|{\cal P}_b(B)|}p_c^{|{\cal P}_c(B)|}.\]

$(ii)$ Let $B\in \bar{{\cal A}^\star(S,G)}$ be a finite DA with over-source the free set $S$. We have
\[`P_{p_a,p_b,p_c}({\A}_\star^{S}=B)=p_a^{|B|}p_b^{|\bar{{\cal P}_{S,b}}(B)|}p_c^{|\bar{{\cal P}_{S,c}}(B)|},\]
where for $d\in\{b,c\}$, $\bar{{\cal P}_{S,d}}(B)={\cal P}_d(A)\cup\{x, x\in S, C_x^\star=d\}$.
\end{pro}
Our gas model of type 2 on $G$ is built as follows.
For any $x\in V$ and $`o\in \Omega$, set
\ben
\label{yop}
X^\star_x(`o)&:=&
\left\{
\begin{array}{ll}
0 & \textrm{ if } C_x(`o)=b\\
1 & \textrm{ if } C_x(`o)=c\\
\dis \min_{c : \textrm{ children of }x} X^\star_{d}(`o)& \textrm{ if }C_x(`o)=a
\end{array}\right.\\
&=&\1_{C_x^\star=c}+ \1_{C_x^\star=a} \min_{d : \textrm{ children of }x} X_{d}(`o).
\een

Once again this recursive definition allows to eventually compute $X^\star_x$ for any $x$, if $p_a<p_{crit}^{\{x\}}$. 
By a simple recursion one easily notices that the values $X^\star_x$ all belong to $\{0,1\}$. 
\begin{remi} Since the values $X^\star_x$ all belong to $\{0,1\}$, the $\min$ operator coincides with the product, and also with the ``and'' operator, interpreting the gas occupation as Boolean variables.
Bousquet-Mélou \cite{MI1} when she considers the case $(p_1,p_2,p_2,p_2)$ in the square lattice case works under an equivalent model even if she uses a different vocabulary (we refer to Section \ref{SQAM} for some hints). 
\end{remi}

\begin{theo}\label{gfdg} Let $G=(V,E)$ be an agreeable graph.\\
$(i)$ For any $x\in V$ and  $p_a<p_{crit}^{\{x\}}$ we have
\[`E_{p_a,p_b,p_c}(X^\star_x)=p_c+`P_{p_a,p_b,p_c}(|{\cal P}_b({\A}^\star)|=0)=p_c+{\bf G}^{x}(p_a,p_c),\]
where ${\bf G}^{x}(u,v)=\sum_{A, |A|>0, {\cal S}(A)=x} u^{|A|}v^{|{\cal P}(A)|}$ is the GF of the set of DA with source $\{x\}$ counted according to their area and perimeter. \\
$(ii)$ Let $S$ be a free finite subset of $G$, and $p_a<p_{crit}^S$. Under $`P_{p_a,p_b,p_c}$, we have
\begin{equation}\label{fondez}
`E_{p_a,p_b,p_c}\l(\prod_{x\in S} X_x^\star\r)=`P_{p_a,p_b,p_c}(X_x^\star=1,x\in S)= \overline{{\bf G}}^{S}(p_a,p_c).
\end{equation}
where $\bar{{\bf G}}^{S}(u,v)=\sum_{A, {\cal S'}(A)\subset S} u^{|A|}v^{|\bar{{\cal P}}(A)|}$  is the GF of the set of DA with over source $S$ counted according to their area and perimeter. 
\end{theo}
\Proof $(i)$ The first equality follows the fact that $X^\star_x=1$ if and only if $C_x^\star=c$ or ${\A}_\star^{\{x\}}$ contains $x$ and has only $c$-colored perimeter sites. For the second equality, write
\be
`P_{p_a,p_b,p_c}(|{\cal P}_b({\A}^{\{x\}}_\star)|=0)&=&\sum_{A, |A|>0, {\cal S}(A)=x} p_a^{|A|}p_c^{|{\cal P}(A)|}= {\bf G}^{x}(p_a,p_c).
\ee 
The proof of $(ii)$ follows the same lines~: if all the values $X_x^\star$ for $x\in S$ equals 1, then the DA $A^{S}$ has all its perimeter sites $c$-colored, and the points of $S\setminus S_{\bullet}$ must be $c$-colored.~$\Box$

\subsection{Comments}

The gas model of type 2 on the square lattice (or on an other lattice), under $`P_p$ for $p<p_{crit}$ is well defined, and one may again study the process $X_x$ when $x$ traverses a choosing line. We can establish a formula similar to \eref{rezr} and a Proposition similar to Proposition \ref{uni} for the gas model of type 2, using the same arguments. In order to compute the density of the gas model of type 2 (this should provide a way to compute $p_{crit}$ which is unknown) a result analogous to Theorem \ref{zoa} is needed: a description of the gas model on a line, or on ``something''. But, the gas model of type 2 even if very similar to that of type 1 has a very different behavior. Except for very special values of $p_a,p_b$ and $p_c$, it is not Markovian on horizontal line (and no Markovian of order 2, and we think no Markovian of order $k$ for any $k$).

We want here to point out that some gas models generalizing the gas model of type 2 may have the wanted property to be Markovian (or have some suitable structural properties). 
An  idea would be to enrich the gas in building a model of gas taking its values in a set larger than $\{0,1\}$ (the gas model of type 2 will appear as a kind of projection of the generalized model), and again use the min operator inside the DA of calculus.  The minimum of some random variables with a prescribed law $\nu$ would replace the minimum of Bernoulli random variables. By chance, maybe there exists a parameter $\nu$ (or a family of distribution $\nu$) for which the gas has a simple description.

\section{Appendix}
\subsection{Proof of Proposition \ref{yopop}} 
Consider a random DA ${\A}^{S}$. Denote by $L_0=S\cap {\A}^{S}$. Then for any $i\geq 1$, set
\[L_i=\{v, v\in {\A}^{S}, v \mbox{ has a father in } L_{i-1}, v\notin \cup_{j\leq i-1}L_j\}.\] 
The random sequence $(L_i)$ gives a decomposition of  ${\A}^{S}$ into layers~: the cells in $L_i$ are the children of the nodes in $L_{i-1}$ that are not children of any cells belonging to $L_j$, $j\leq i-2$.
We want to prove that if $p<1/K$  then $`P_p$ a.s., there exists $k$ such that $|L_k|=0$. \par
The number $N_c$ of children in $A$ of a given cell $c$ of $A$ has law Binomial$(o(c),p)$, where $o(c)$ is the out degree of $c$ in $G$, that is
\[`P_p(N_c=j)=\binom{o(c)}j p^j(1-p)^{o(c)-j} \textrm{ for }j \in \{0,\dots, o(c)\}.\]
Therefore, given $L_i=l_i$ for $i\leq k$, the random variable $|L_{k+1}|$ has law Binomial$(o(C),p)$ where $o(C)$ is the number of children of the nodes of $l_k$ in $G$ that are not children of any nodes of the $l_i$,$i<k$. 
Hence, given $L_i=l_i$ for $i\leq k$, the random variable $|L_{k+1}|$ is smaller for the stochastic order than Binomial$(K\, |l_k|,p)$~: we say that $X$ is smaller than $Y$ for the stochastic order if $`P(X\geq x)\leq `P(Y\geq x)$ for any $x$ (we write $X\leq _S Y$). This is simple property of the binomial distribution 
clear via the representation of binomial random variables as sum of i.i.d. Bernoulli random variables.
A simple iteration, shows that $|L_k|$ is smaller for the stochastic order than $Z_k$ where $Z_0=|L_0|$, and where given $Z_{k-1}=z_{k-1}$, $Z_k$ is Binomial$(Kz_{k-1},p)$ distributed. There exists a probability space, on which one may construct two sequences $(l'_k)$ and $(Z'_k)$ such that 
\[(l'_k)\sur{=}{(d)}(|L_k|) \textrm{ and } (Z_k')\sur{=}{(d)}(Z_k)\]
and such that a.s. $l'_k\leq Z'_k$. On this space, $\inf\{k, l'_k=0\}\leq T:=\inf\{k, Z'_k=0\}$.\par 
But, $Z_k'$ is simply a Galton-Watson process with offspring distribution Binomial$(K,p)$ (see Athreya \& Ney \cite{AtN}), and then it is classical that such a Galton-Watson process eventually dies out ($T<+\infty$ a.s.) if the mean of its offspring distribution is smaller than 1 (if $|L_0|$ is finite, which is the case here); here it is the case if $p\leq 1/K$.  
~$\Box$ 
\begin{com}\label{soupe}If $G$ is a tree in which at depth $k$, all the nodes have $k+1$ children, then $p_{crit}=0$. Indeed, under $`P_p$, the mean number of children of an individual having depth larger than $1/p$ is larger than 1. In the following levels, the number of individuals are larger for the stochastic order than a super-critical Galton-Watson process, for which $T=+\infty$ with positive probability.
\end{com}

\subsection{Square lattice DA counted according to the area: The historical approach}
\label{SQAM}

The content of this section \ref{SQAM} is a reformulation of a part of Bousquet-Mélou \cite{MI1} which is in some sense the mathematical justification of the work of Dhar \cite{DH6} and concerns the enumeration of DA on the square lattice using a gas model of type 1 (this terminology is not introduced there).

First, instead of working on the square lattice, a cyclic directed square lattice having $n$ cells in each row is introduced~:
 the cylinder $\Cy(n) = (\mathbb{Z}/n\mathbb{Z})\times \mathbb{Z}$, for $n\geq 2$ which is the lattice $\Sq$ quotiented by a congruence relation~:  $(x,y)$ has children $(x \mod n ,y+1)$ and $(x+1 \mod n ,y+1)$.

Consider a DA $A$ with source $C$ in $\Cy(n)$. 
\begin{figure}[htbp]
\psfrag{v}{$\varnothing$}\psfrag{gi}{$G_i$}\psfrag{fi}{$G_{i+1}$}\psfrag{1}{1}\psfrag{2}{2}\psfrag{3}{3}\psfrag{31}{31}\psfrag{i}{row $i$}\psfrag{i+1}{row $i+1$}\psfrag{12}{12}
\centerline{
\includegraphics[height=1.8cm]{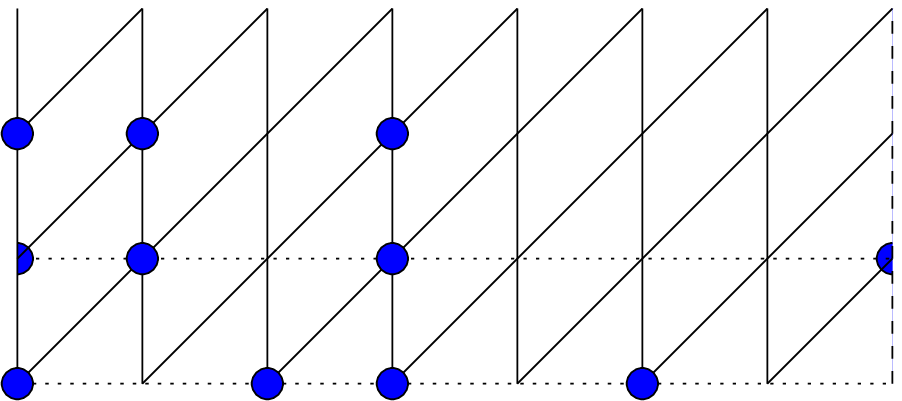}
}
\captionn{\label{cy}A DA on the cylinder. The maximum source has 4 cells, the area is 10. The right most and left most vertices are identified.} 
\end{figure} 
Removing the first row of $A$, the following formula is easily checked,
\begin{equation}
\label{eq1}
{\G_{C}^{{\sf cy}}}(x)=x^{|C|}\sum_{D \subset {\cal N}(C)} \G_{D}^{{\sf cy}}(x),
\end{equation}
where ${\cal N}(C)=\{i,i+1 \mod n,i\in C\}$ is the perimeter sites of $C$.\par

Then, Bousquet-Mélou introduced a gas model equivalent to the gas model of type 1. The rows of the cylinder are labeled by integers, the $i+1$ row being above the $i$th (see Fig. \ref{cy}). 
The gas occupation on a row is random vector $X=(X_1,\dots,X_n)$ taking its values in $\{0,1\}^n$; its distribution is described by
$F_C=`P(X_i=1,i\in C, X_i=0, i\in {}^c C)$ or by $f_C=`P(X_i=1,i\in C)$ using the inclusion-exclusion principle. 
For any $j$, the vector $X^j$ giving the distribution of the gas on the $j$th row is given from the $X^{j+1}$th by the following stochastic evolution (gas model of type 1): \par
\[X^j_i=B_i^j(p) (1-X^{j+1}_i)(1-X^{j+1}_{(i+1 \mod n)})~~~\textrm{ for any }i\in \mathbb{Z}/n\mathbb{Z} ,j\in\mathbb{Z}\]
where the $B_i^j(p)$ are i.i.d. Bernoulli$(p)$ random variables (in other words, a cell $i$ in the $j$th row is occupied with probability $p$ if and only if $i$ and $i+1$ are empty in the $j+1$th cell, and if a Bernoulli$(p)$ tossing is a success, the Bernoulli tossing being independent).

One gets $f^j_C(p)=p^{|C|}`P(X^{j+1}_i=0, i\in {\cal N}(C))$ and
by the inclusion-exclusion principle,
\begin{equation}
f^{i}_C(p)=p^{|C|}\sum_{D \subset {\cal N}(C)}(-1)^{|D|} f^{i+1}_D(p),
\label{eq2}
\end{equation}
this in terms of $F^{(i)}$ reads
\[F^i_C(p)=\l(\frac{p}{1-p}\r)^{|C|}\sum_{D\in \{1,\dots,n\}\setminus {\cal N}(C)} (1-p)^{n-|{\cal N}(D)|}
F^{i+1}_D(p).\]
The sequence $(X^j)$ forms a Markov chain which is aperiodic, irreducible with a finite number of states. Hence, there exists a unique invariant distribution $F$;  the corresponding function $f$ satisfies
\begin{equation}
f_C(p)=p^{|C|}\sum_{D \subset {\cal N}(C)}(-1)^{|D|} f_D(p).
\label{eq3}
\end{equation}
Now, the system \eref{eq3} is solved explicitly by checking that 
\begin{equation}\label{eer}
F_D(p)=\frac{1}{Z_n}\l(\frac{p}{1-p}\r)^{|D|}(1-p)^{|{\cal N}(D)|} \textrm{ where }Z_n=\sum_{D\subset\{1,\dots,n\}}\l(\frac{p}{1-p}\r)^{|D|}(1-p)^{|{\cal N}(D)|}
\end{equation}
is solution. Notice that for each $C$, $F_C$ is a rational function of $p$ since it is solution of a linear system with polynomial coefficient in $p$.
Observe \eref{eq1} and \eref{eq3}. If $f$ is a solution of \eref{eq3} then 
\[\G_{C}^{{\sf cy}}(-p) = (-1)^{|C|}f_C(p)\]
is a solution of \eref{eq1} (in fact there is only one solution up to a multiplicative constant, which is fixed taking $-\sum_C \G_{C}^{{\sf cy}}(-p)=1$).
Denote by $\G_{c}^{{\sf cy}}$ the series of DA with source $c$, a unique cell. We have 
\[-\G_{c}^{{\sf cy}}(-p)=f_{c}(p)=`P(X_1=1),\] 
the so-called \it density \rm of the gas model.
The series $\G_{c}^{{\sf cy}}$, we should write $\G_{c}^{{\sf cy}(n)}$, is the series of DA in the cylinder. It is clear that its first $n$ coefficients coincide with that of the series $\G_c^{{\sf sq}}$ of DA on the square lattice with source $O$. Hence in the space of formal series
\[\G_{c}^{{\sf cy}(n)} \sous{\tend}{n\to \infty} \G_{c}^{\Sq}.\]
The explicit solution of the gas model on $\Cy(n)$ is derived as follows~: using that $|{\cal N}(D)|=|D|+|{\cal N}_r(D)|$ where ${\cal N}_r(D)=\{i, i\in D,i+1\notin D\}$, $Z_n$ reads $\sum_{D\subset\{1,\dots,n\}}{p}^{|D|}(1-p)^{|{\cal N}_r(D)|}$ and  is guessed (or viewed) to have the following form
\[Z_n=\sum_{(y_1,\dots,y_n)\in\{0,1\}^n}\prod_{i=1}^n V(y_i,y_{i+1}) \textrm{ where }\]
$V(0,0)=V(0,1)=1$, $V(1,0)=p(1-p)$, $V(1,1)=p$ and $y_{n+1}:=y_1$. 
This corresponds in some sense to a ``cyclical Markovian model'' where a transition $0\to1$ or $0\to0$ is ``counted'' $1$, a transition $1\to0$ is counted $p(1-p)$, a transition $1\to1$, counted $p$. 
In other words, since $y_{n+1}=y_1$, $Z_n$ is  $tr(Y^n)$ where  
$Y=\left(
\begin{array}{ll}
1&1\\p(1-p)&p
\end{array}
\right);$ 
this is equal to $\lambda_1^n+\lambda_2^n$ where $\lambda_1>\lambda_2$ are the two eigenvalues of $Y$, that are~: $\lambda_{1,2}=\frac{1+p\pm\sqrt{1+2p-3p^2}}2$. 
It remains to express the density. In Bousquet-Mélou \cite{MI1} this is done by differentiation of \eref{eer} (where some mute variables replace $(p/(1-p))$ and $(1-p)$). One may also compute this density by computing 
$W_n/Z_n$ where $W_n= (Y^n)_{2,2}$  which corresponds to the contribution of the configurations where the first cell is occupied (in other words, $y_1=y_{n+1}=1$). With very simple linear algebra, we get $W_n=(\lambda_1^n(\lambda_1-1)+\lambda_2^n(\lambda_2-1))/(\lambda_2-\lambda_1)$. 
Since $\lambda_1>\lambda_2$, 
\begin{equation}\label{more}
W_n/Z_n\to (1-\lambda_1)/(\lambda_2-\lambda_1)
\end{equation}
which is equal to $\frac{1}{\alpha_\bullet}/(\frac{1}{\alpha_\bullet}+\frac{1}{\alpha_\circ})$.  The convergence in \eref{more} is considered in the formal sense.

\begin{com} In fact, the solution of \eref{eer} is more or less guessed by Dhar \cite{DH5} who works on the whole lattice. He ``notices'' that the distribution of particles on two consecutive lines coincides with the hard particle distribution of activity $p/(1-p)$ on two consecutive lines. On the cylinder, this can be interpreted saying that the probability to see $j$ occupied cells on $G_i\cup G_{i+1}$ is $\lambda (p/(1-p))^j$ where $\lambda$ is a constant. Given this, the distribution on a line is then just a marginal of this distribution, easy to compute.
\end{com}

In \cite{MI1}, Bousquet-Mélou generalized the study of gas model allowing some new evolutions between lines. With these tools, she is able to give some formal link  between density of a gas and GF of DA on the square lattice (and also on other lattices) counting DA according to several parameters, the area, the perimeter, the right perimeter and the ``loops''. These links are once again formal. 
For example, she derived the following formula~:
the area and perimeter GF for one-source DA on the cyclic square lattice is 
\begin{equation} \label{ser}
\sum_{A}t^{|A|}x^{|{\cal P}(A)|}=1-x-\rho(p_1,p_2,p_2,p_2)
\end{equation}
where $p_1=1-x-t$ and $p_2=1-x$. The quantity $\rho(p_1,p_2,p_3,p_4)$ is the density of a model of gas  obtained by the rules evolution given in the following table. 
\[\begin{array}{|c|c|c|c|c|}
\hline
\textrm{Occupation of $x$ and $x+1$ in $Z_1$} & 0 \textrm{ and } 0 &  1 \textrm{ and } 0 & 0 \textrm{ and } 1 &1 \textrm{ and } 1\\
\hline
\textrm{Occupation of $x$ in $Z_0$} & \textrm{Bernoulli}(p_1)  & \textrm{Bernoulli}(p_2) & \textrm{Bernoulli}(p_3)  & \textrm{Bernoulli}(p_4).  \\ \hline
  \end{array}\]

\small 
\bibliographystyle{plain}

\end{document}